\documentclass[a4paper,11pt]{book}
\usepackage[centertags]{amsmath}
\usepackage{amscd}
\usepackage{amsthm}
\usepackage{amssymb}
\usepackage{epigraph} 
\usepackage{enumerate}
\usepackage{multicol}
\usepackage{esvect}
\usepackage[absolute,overlay]{textpos}
\usepackage{bm}

\usepackage[english,catalan]{babel}
\usepackage{color}
\usepackage{tikz}
\usepackage{tabls} %per tenir m{\'e}s control sobre l'espai en les taules nom{\'e}s cal aquest paquet
\usepackage{indentfirst} %per comen\c{c}ar primer parragraf de seccio o capitol amb indent
\usepackage{mathtools}

\DeclarePairedDelimiter\floor{\lfloor}{\rfloor}

\usepackage[T1]{fontenc}
\usepackage[sc]{mathpazo}
\usepackage[latin1]{inputenc}
\usepackage[twoside,bindingoffset=1cm]{geometry}
\usepackage{array,amsfonts}

\newtheorem{defn0}{Definition}[chapter]
\newtheorem{prop0}[defn0]{Proposition}
\newtheorem{thm0}[defn0]{Theorem}
\newtheorem{lemma0}[defn0]{Lemma}
\newtheorem{corollary0}[defn0]{Corollary}
\newtheorem{example0}[defn0]{Example}
\newtheorem{remark0}[defn0]{Remark}
\newtheorem{conjecture0}[defn0]{Conjecture}

\newenvironment{definition}{ \begin{defn0}}{\end{defn0}}

\newenvironment{theorem}{\bigskip \begin{thm0}}{\end{thm0}}
\newenvironment{lemma}{\bigskip \begin{lemma0}}{\end{lemma0}}
\newenvironment{corollary}{\bigskip \begin{corollary0}}{\end{corollary0}}

\usepackage{fancyhdr}
\newif\ifprivate
\privatetrue

\def\???{\ifprivate {\bf {???}} \marginpar{\begin{center}{\Huge {\bf ?}}\end{center}}
\else \fi}
\begin{document}

\pagestyle{empty}

\begin{titlepage}
\begin{center}
\begin{figure}[htb]
\begin{center}
\includegraphics[width=6cm]{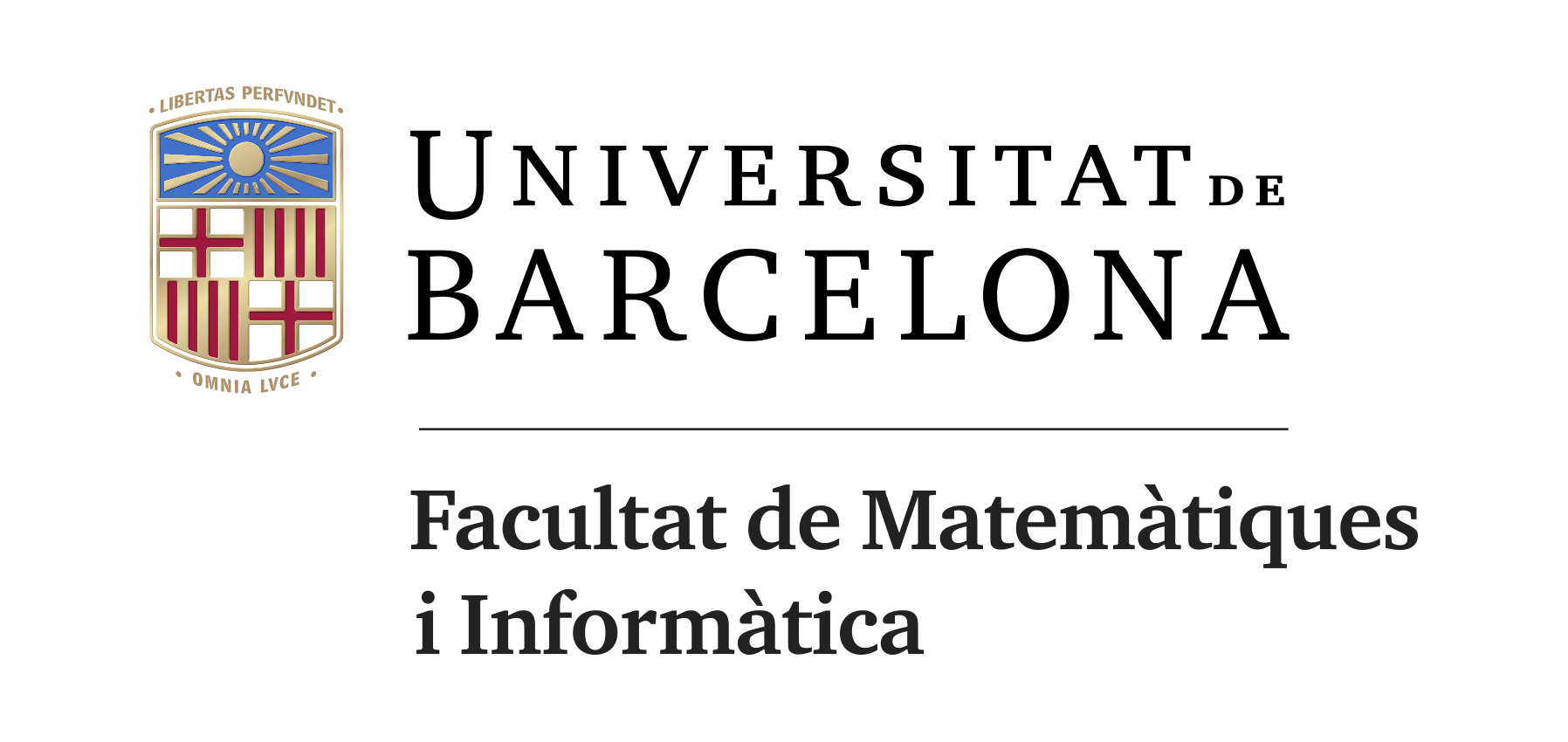}
\end{center}
\end{figure}

\vspace*{1cm}
\textbf{\LARGE GRAU DE MATEM\`{A}TIQUES } \\
\vspace*{.5cm}
\textbf{\LARGE Treball final de grau} \\

\vspace*{1.5cm}
\rule{16cm}{0.1mm}\\
\begin{Huge}
 \textbf{The real numbers\\ in inner models \\ of set theory} \\
\end{Huge}
\rule{16cm}{0.1mm}\\

\vspace{1cm}

\begin{flushright}
\textbf{\LARGE Autor: Mart\'{i}n Soto Quintanilla}

\vspace*{2cm}

\renewcommand{\arraystretch}{1.5}
\begin{tabular}{ll}
\textbf{\Large Director:} & \textbf{\Large Dr. Joan Bagaria Pigrau } \\
\textbf{\Large Realitzat a:} & \textbf{\Large  Departament de Matem\`{a}tiques i Inform\`{a}tica   } \\
 & \textbf{\Large } \\
\\
\textbf{\Large Barcelona,} & \textbf{\Large 13 de juny de 2022 }
\end{tabular}

\end{flushright}

\end{center}

\end{titlepage}

%%%%%%%%%%%%%%%%%%%%%%%%%%%%%%%%%%%%%%%%%%%%%%%%%%%%%%%%%%%%%%%%%%%%%%%%%

\pagestyle{fancy}
\newpage \thispagestyle{empty}

\section*{Abstract}
We study the structural regularities and irregularities of the reals in inner models of set theory. Starting with $L$, G\"{o}del's constructible universe, our study of the reals is thus two-fold. On the one hand, we study how their generation process is linked to the properties of $L$ and its levels, mainly referring to \cite{Gaps}. We provide detailed proofs for the results of that paper, generalize them in some directions hinted at by the authors, and present a generalization of our own by introducing the concept of an infinite order gap, which is natural and yields some new insights. On the other hand, we present and prove some well-known results that build pathological sets of reals.

We generalize this study to $L[\#_1]$ (the smallest inner model closed under the sharp operation for reals) and $L[\#]$ (the smallest inner model closed under all sharps), for which we provide some introduction and basic facts which are not easily available in the literature. We also discuss some relevant modern results for bigger inner models.

{\let\thefootnote\relax\footnote{2020 Mathematics Subject Classification. 03E45, 03E15, 03C62}}

\newpage
\thispagestyle{empty}

\section*{Agra\"{i}ments}

Estic immensament agra\"it al Dr. Joan Bagaria per tota la seva ajuda, atenci\'{o} i paci\`{e}ncia. Aprendre d'ell ha estat un gran plaer, i una font de motivaci\'{o} constant. Gr\`{a}cies als seus coneixements i passi\'{o} que altruistament comparteix, la Teoria de Conjunts \'{e}s per mi un m\'{o}n m\'{e}s ric i profund que mai.

Voldria tamb\'{e} dedicar aquest treball a totes les meves companyes del doble grau. Han fet d'aquests cinc anys els millors de la meva vida, i sense elles no seria on s\'{o}c ni qui s\'{o}c. Cada classe, cada conversa i cada mirada ho han estat tot per mi.

I gr\`{a}cies infinites als meus pares per haver-m'ho donat tot, i per haver-me ensenyat a aprendre i a ser persona. 

\newcommand{\chapquote}[3]{\begin{quotation} \textnormal{#1} \end{quotation} \vspace{-0,8cm} \begin{flushright} -- #2, \textit{#3}\end{flushright} }

\selectlanguage{english}
\pagenumbering{roman} \setcounter{page}{0}
\tableofcontents
 \thispagestyle{empty}
 \newpage
  \thispagestyle{empty}
  \begin{center}
      \scalebox{1.2}{
      \begin{tikzpicture}
        \def\nangle{70} \def\r{7.2} \def\mygreen{green!50!black}
        
        \foreach \i/\lbl in {
            2/1,
            4/2,
            6/3,
            8/\omega,
            9/\omega+1,
            10/\omega+2,
            14/2\omega,
            16/3\omega,
            20/\omega\times\omega}          
            {
            
            }
            
            \path (0,0) -- (0,7) node[pos=-0.05]
            {$0=\varnothing$}
            node[pos = 0.12, yshift=0pt,fill=black,circle,minimum size=3pt, inner sep=0pt] {}
            node[pos = 0.18, yshift=0pt,fill=black,circle,minimum size=3pt, inner sep=0pt] {} 
            node[pos = 0.24, yshift=0pt,fill=black,circle,minimum size=3pt, inner sep=0pt] {} 
            node[pos=0.32,sloped] {$\dots$}
            node[pos = 0.4, yshift=0pt,fill=black,circle,minimum size=3pt, inner sep=0pt] {} 
            node[pos = 0.46, yshift=0pt,fill=black,circle,minimum size=3pt, inner sep=0pt] {} 
            node[pos = 0.52, yshift=0pt,fill=black,circle,minimum size=3pt, inner sep=0pt] {} 
            node[pos=0.6,sloped] {$\dots$}
            node[pos = 0.68, yshift=0pt,fill=black,circle,minimum size=3pt, inner sep=0pt] {} 
            node[pos=0.76,sloped] {$\dots$}
            node[pos = 0.84, yshift=0pt,fill=black,circle,minimum size=3pt, inner sep=0pt] {} 
            node[pos=0.92,sloped] {$\dots$}
            ;
            
            \path (0.3,0) -- (0.3,7) node[pos=0.12]
            {1} node[pos=0.18]
            {2} node[pos=0.24]
            {3} node[pos=0.4, xshift = 0.46cm]
            {$\omega = \mathbb{N}$} node[pos=0.46, xshift= 0.35cm]
            {$\omega+1$} node[pos=0.52, xshift= 0.35cm]
            {$\omega+2$} node[pos=0.68, xshift= 1.1cm]
            {$\omega+\omega = 2\times\omega$} node[pos=0.84, xshift= 0.4cm]
            {$\omega\times\omega$};
            
            \draw (-2.2, 5) node{$\mathbb{Q}$};
            \draw (-1.6, 5.8) node{$\{\pi, \frac{3}{5}\}$};
            \draw (-2.9, 5.6) node{$\mathbb{R}$};
        
        \draw[line width=2pt, cap=butt, rounded corners] (90+0.5*\nangle:1.1*\r) -- (0,0) node[yshift=1pt,fill=black,circle,minimum size=5pt, inner sep=0pt] {} -- (90-0.5*\nangle:1.1*\r);
    \end{tikzpicture}
    } \\
     \textit{A usual representation of the set theoretic universe, with the \textnormal{ordinals} as its backbone.}\vspace{1.2cm}
    \scalebox{1.2}{
  	\begin{tikzpicture}{
        \def\nangle{70} \def\r{7.2} \def\mygreen{green!50!black}
        
        \foreach \i/\lbl in {
            2/1,
            4/2,
            6/3,
            8/{\omega},
            9/{\omega+1},
            10/{\omega+2},
            14/{2\omega},
            16/{3\omega},
            20/{\omega\times\omega}}          
            {
            
            }
            
            \path (0,0) -- (0,7) node[pos=-0.05]
            {$0=\varnothing$}
            node[pos = 0.14, yshift=0pt,fill=black,circle,minimum size=3pt, inner sep=0pt] {}
            node[pos=0.19, xshift=0.1cm]
            {$\omega=\aleph_0$}
            node[pos=0.32,sloped] {$\dots$}
            node[pos = 0.4, yshift=0pt,fill=black,circle,minimum size=3pt, inner sep=0pt] {}
            node[pos=0.5,sloped] {$\dots$}
            node[pos = 0.58, yshift=0pt,fill=black,circle,minimum size=3pt, inner sep=0pt] {}
            node[pos=0.68,sloped] {$\dots$}
            node[pos = 0.76, yshift=0pt,fill=black,circle,minimum size=3pt, inner sep=0pt] {}
            node[pos=0.84,sloped] {$\dots$};
            
            \path (0.3,0) -- (0.3,7)
            node[pos=0.4, xshift = 0.55cm]
            {$\omega_1 = \aleph_1$} node[pos=0.58, xshift= 0.1cm]
            {$\aleph_2$} node[pos=0.76, xshift= 0.1cm]
            {$\aleph_\omega$};

        \draw[line width=1pt] (90+0.5*\nangle:0.3*\r) --  (90-0.5*\nangle:0.3*\r);
        
        \draw[line width=2pt, cap=butt, rounded corners] (90+0.5*\nangle:1.1*\r) -- (0,0) node[yshift=1pt,fill=black,circle,minimum size=5pt, inner sep=0pt] {} -- (90-0.5*\nangle:1.1*\r);}
    \end{tikzpicture}
    } \\
    \begin{textblock*}{5cm}(12cm,22cm) % {block width} (coords) 
   \Large $ \omega_\alpha = \aleph_\alpha$
\end{textblock*}
\vspace{0.2cm}
    \textit{Zooming out beyond the countable ordinals, we see more infinite \textnormal{cardinals}.}
    \end{center}
%%%%%%%%%%%%%%%%%%%%%%%%%%%%%%%%%%%%%%%%%%%%%%%%%%%%%%%%%%%%%%%%%%%%%%%%%

%%%%%%%%%%%%%%%%%%%%% aix{\`o} pels headings %%%%%%%%%%%%%%%%%%%%%%%%

%%%%%%%%%%%%%%%%%%%%%%%%%%%%%%%%%%%%%%%%%%%%%%%%%%%%%%%%%%%%%%%%%%%%%%%%%
\chapter*{Introduction}
\markboth{Introduction}{Introduction}
\addcontentsline{toc}{chapter}{Introduction}

\chapquote{"Sometimes it seems that there is such a complete lack of rule-governed behavior that [some phenomena] just \textit{aren't} rule-governed. But this is an illusion---a little like thinking that crystals and metals emerge from rigid underlying laws, but that fluids or flowers don't."}{Douglas R. Hofstadter}{G\"{o}del, Escher, Bach}

\section*{Real numbers}

The set of real numbers is undoubtedly one of the most important objects in mathematics. The real numbers capture our understanding of continuous space, and so play a crucial role in applied mathematics or mathematical branches like Topology and Geometry. But even in Set Theory, the standard foundation for all of mathematics, this set was central in shaping the developments of the field, and has been thoroughly established as the kernel of many still open questions regarding the nature of the mathematical universe.

\vspace{0,5cm}

The reals played a historically privileged role in the discretization of mathematics. As Kanamori writes \cite{Kan}, up until the ending of the 19$^{\textnormal{th}}$ century mathematics had enjoyed intensional interpretations and methods, where the focus was to be found in processes and empirical intuitions rather than objects (similarly to how physical phenomena or natural intuitive reasoning are most easily understood). And so for instance functions where understood as rules rather than collections, and many objects defined by their properties rather than a construction.

The increasingly proof theoretic approach striving for rigor brought to bear the shortcomings of this approach. Most prominently, the study of limits and continuity made apparent the need for an extensional understanding of functions as acting on points. Collections of discrete objects were especially useful as an underpinning for these previously non-rigorous concepts, such as the well-known constructions of the continuum by Dedekind or Cantor. Of course, Cantor would go on to establish the foundations of Set Theory, and his methods would generalize to the extensional understanding of all mathematical objects as sets we use today.

But the continuum didn't only play an important role because the principal pressures towards an extensional approach came from real analysis. It also did so because of the consequences its formulation revealed. Cantor's discovery of the (infinitely many) different sizes of infinity was met with skepticism. This was the first instance of a common phenomenon: as mathematically useful as the extensional underpinning might be, some of its discrete properties seem counter to our intensional intuitions about how an object should behave. 

Of course, the different sizes of infinity have come to be widely accepted as witnessing the richness of the mathematical universe. But as we see in Chapter 2, other apparent implausibilities, more concrete and pathological, still perplex us. We also notice that this tension between a concept and its underpinning is especially noticeable in the continuum because it is our paradigm example for continuity, and thus its discretization strikes us as less plausible.

\vspace{0,5cm}

When regarding Set Theory not just as the foundation for all of mathematics, but as the self-contained mathematical discipline studying the infinite, the set of reals also plays a crucial role: settling its cardinality has been one of the driving problems of the field since its beginnings. Since it is independent of $ZFC$, the answer can only be motivated by the practical or intuitive value of adding certain axioms, and this is the subject of much of modern Set Theory through Large Cardinal axioms. But how can it be that such a fundamental question is not settled by this otherwise fruitful Set Theory?

To answer that, let us first note that in Logic and Set Theory we identify $\mathbb{R}$ with $\mathcal{P}(\omega)$, and work with the latter. This is not only because $\mathbb{R}$ will always have that cardinality, but also because, regardless of the construction of $\mathbb{R}$ we choose, there will be a very natural correspondence between its members and the subsets of $\omega$, and thus all results will be directly translatable (for instance, we can consider the characteristic function of a subset of $\omega$ as an infinite sequence of binary digits, the binary expression of a real number).

The reason why $ZFC$ doesn't decide the Continuum Hypothesis is that it generally provides very little information about how the Power Set operation behaves. The Power Set axiom tells us we can always apply it, but gives us no information about its richness or its properties, and so other than Cantor's Theorem ($|\mathcal{P}(x)| > |x|$) and some results regarding cofinalities, cardinal exponentiation remains mainly undecided. We need to pin down how the power set operation actually looks to get more information, and inner models do so.

\section*{Inner models}

Inner models are transitive set theoretic classes that satisfy the axioms of $ZF$ and contain all the ordinals. That is, a smaller mathematical universe that nonetheless contains the same ordinal backbone. \\
\begin{center}
\vspace{-4.5cm}
		\begin{tikzpicture}
        \def\nangle{80} \def\sangle{40} \def\r{3} \def\mygreen{green!50!black}
            
            \path (0,0) -- (0,7) 
            node[pos = 0.05, yshift=0pt,fill=black,circle,minimum size=3pt, inner sep=0pt] {}
            node[pos=0.12,sloped] {$\dots$}
            node[pos = 0.18, yshift=0pt,fill=black,circle,minimum size=3pt, inner sep=0pt] {}
            node[pos=0.24,sloped] {$\dots$}
            node[pos = 0.3, yshift=0pt,fill=black,circle,minimum size=3pt, inner sep=0pt] {}
            node[pos=0.36,sloped] {$\dots$}
            ;

        \draw (1.5, 2.95) node {$M$};
        \draw (-0.5, 2.1) node {$\mathbb{R}^M$};
        \draw (-1.6, 2.07) node {$\mathbb{R}$};
        \draw (3, 2.95) node {$V$};
        
        \draw[line width=2pt, cap=butt, rounded corners] (-3,2.7) -- (0,0) node[yshift=1pt,fill=black,circle,minimum size=5pt, inner sep=0pt] {} -- (3, 2.7);
        \draw[line width=2pt, cap=butt, rounded corners] (-1.5, 2.7) -- (0,0) node[yshift=1pt,fill=black,circle,minimum size=5pt, inner sep=0pt] {} -- (1.5, 2.7);
    \end{tikzpicture}\\
    \textit{Inner model M (with set of reals $\mathbb{R}^M$) built inside universe V}
    \vspace{0.1cm}
\end{center}
Each one of these models will have certain concrete properties (that is, satisfy some axioms additional to $ZF$), and so decide way more statements. In fact, as we'll see, inner models usually present some kind of structural regularity that facilitates their study.

This regularity and the minimality of the model is in great part achieved by concretely defining what the power set function does (especially in canonical inner models, see Chapter 3). For instance, G\"{o}del's $L$, the smallest inner model, can be understood as implementing the idea of choosing the simplest possible power set function allowed by the expressive power of our language. That is, at every level of its hierarchy, we add only what we can talk about, and thus must be present for the model to be coherent. Other bigger models implement in their construction slightly less simple power set operations.

So it is natural that in these models the reals see their complexity greatly reduced. The issue is, this yields some unwanted behaviour that crashes with our intensional understanding of the reals. That is, as useful as the regularity properties are for some set theoretic purposes, these models can't wholly capture an inherently complex object like the reals.

We will study how the regularity of a model affects its reals, their resulting unwanted properties, and the ways in which they are correlated.

\section*{This work}

Our study of the reals of a model thus takes two directions. On the one hand, we study how the reals of the model are generated, and how this generation process is intertwined with the construction of the model itself. On the other hand, we study some apparently pathological properties of sets of reals inside that model.

In Chapter 1 we develop the first direction for $L$. We mainly study the generation process through the concepts and developments of \cite{Gaps}, while still proving other results and tracing some necessary tools to other authors. We provide detailed proofs for the main theorems of the paper, generalize the results in certain directions hinted at by the authors and present some generalizations of our own. Most prominently, we introduce the notion of an infinite order gap, which helps phrase some natural questions, and whose study yields some further insight into the model-theoretical structure of the hierarchy. So all in all this chapter has a considerable mathematical load. We note nonetheless that our presentation of the paper isn't exhaustive, as some of its concrete results weren't as relevant to our main purpose and thus have been excluded.

In Chapter 2 we develop the second direction for $L$. We present and prove some well-known results, and remark upon the possible heuristic reasons for their apparent plausibility issues.

In Chapter 3 we generalize the previous study to bigger inner models. We present and motivate the models $L[\#_1]$ and $L[\#]$ in some detail, and comment on how the two directions generalize to them. We finally reflect also upon the situation for more complex inner models, by discussing some more general and modern results in Inner Model Theory and Descriptive Set Theory relevant to our purposes.

\vspace{0,5cm}

Since the topics involved are slightly advanced, we can't provide exhaustively all preliminaries needed, so a minimal familiarity with basic concepts and results in set theory and model theory is required to completely understand some of the mathematics (especially in Chapter 1, and some parts of Chapter 3). But we do always state the results used and provide references for them which contain the required background.

\mainmatter
\chapter[Gaps in L]{Gaps in L}

G{\"o}del's constructible universe, $L$, a proper class transitive $\in$-model, was devised to prove some independence results. Indeed, even though neither $CH$ nor $AC$ are needed to construct it, $L$ is a model of $ZFC$ and the Generalized Continuum Hypothesis, so both $AC$ and $GCH$ were proven irrefutable from $ZF$ (assuming of course $ZF$'s consistency). The theory of $L$ is especially tame and simple, and easily resolves many questions independent of $ZF$. $L$ is in fact the smallest inner model of $ZF$, that is, any other inner model contains it. So its canonicity ensured that deeper aspects of its structure would become object of fruitful study. The most celebrated results on this area are Jensen's systematical Fine Structure Theory and the study of $L$-indiscernibles which culminated in the discovery of sharps (used in Chapter 3).

Here we exploit some of these advances, most concretely using the framework and results of \cite{Gaps}, to study the generation process of 	$\mathbb{R}$ in $L$ and understand the role of the constructible reals in the hierarchy levels. The reals turn out to be tightly interlinked with model theoretic properties affecting the whole hierarchy, thanks to their ability to codify set theoretic information. We end up with a picture of the extreme regularity which the Axiom of Constructibility ($V=L$) bestows upon the set theoretic universe, as seen through the reals.

Now, $\omega$ is an object easily definable and of extreme absoluteness: its definition is $\Sigma_0$, and so all transitive models of $ZF$ identify $\omega$ correctly. So when investigating how the members of $\mathcal{P}(\omega)$ behave in $L$ (or any inner model), we are really studying the properties of the power set function in this model (since those of $\omega$ remain unchanged). In fact, as mentioned above, $V=L$ can be understood as artificially determining the power set operation to be the simplest one possible, so it's not surprising that $\mathcal{P}(\omega)$ is so deeply affected by it, or that pursuing this study leads to a generalization to iterations of $\mathcal{P}$ and its values on other cardinals, as presented further down.

\newpage

\section{Preliminaries}

We work in the language of first order logic plus the binary predicate $\in$, all further symbols serving as abbreviations, and assume the axioms of $ZF$. An introduction to $L$, and an exposition of the following well-known results of model theory and set theory, can be found in resources like \cite{Dev} and \cite{Kun}.

\begin{definition}
By $Def^A (P)$ we mean the set of all elements definable over $\langle A, \in \rangle$ with parameters in $P \subseteq A$. That is, all elements $a \in A$ such that, for some formula $\phi$, natural number $k$ and string of parameters $\bar{b} \in [P]^k$, $\{ a \} = \{ x\in A \: | \: \langle A, \in \rangle \models \phi(x, \bar{b}) \}$. \\ We abbreviate $Def^A = Def^A (\varnothing)$. We also say $A$ is pointwise definable from $P\subseteq A$ if $Def^A (P) = A$, and pointwise definable if $Def^A = A$.
\end{definition}

We also abbreviate $\langle A, \in \rangle$ to $A$, since we mostly consider $\in$-models, and it's always clear from context when we use $A$ as a model.
\vspace{-0,2cm}
\begin{theorem}
{\bf (Skolem hull argument) (I.15.28 in  \cite{Kun})}
Assume there is a binary relation $R$ that well-orders $A$ and is definable over $A$ with parameters in $P\subseteq A$. Then $Def^A (P) \prec A$.
\end{theorem}

\begin{theorem}
{\bf (Downward L\"{o}wenheim-Skolem) (I.15.10 in  \cite{Kun})}
Fix an infinite set $A$ and an infinite cardinal $\kappa \leqslant |A|$, and $S \subseteq A$ with $|S| \leqslant \kappa$. Then there's a $B \subseteq A$ with $B \prec A$, $S \subseteq B$ and $|B| = \kappa$.  
\end{theorem}

\begin{theorem}{\bf  (Global well-order in \textit{L}) (II.3.2, 3.3 in  \cite{Dev}, and \cite{Boo})}
\label{bottLine}
There is a formula $W(x, y)$ such that for every ordinal $\alpha$
\begin{itemize}
    \item[] $L_{\alpha+5}$ $\models$ "W(x, y) well-orders $L_\alpha$"
\end{itemize}
Moreover, for every limit ordinal $\alpha$
\begin{itemize}
    \item[] $L_\alpha$ $\models$ "W(x, y) well-orders every $L_\beta$, and thus the universe"
\end{itemize}
Furthermore, for every $n\in\omega$ there is a formula $W_n$ without parameters such that, if $\alpha = \lambda + n$, then
\begin{itemize}
    \item[] $L_{\alpha}$ $\models$ "$W_n$(x, y) well-orders the universe"
\end{itemize}
\end{theorem}
 
 This last improvement on definability (a construction by Boolos) will prove vital for fine-structural purposes. It is basically achieved by using a flat pairing function for the successor case (see Definition 1.17). We call the $L_\alpha$-definable relation represented by this formula $\leqslant_{L_\alpha}$ (and use $<_{L_\alpha}$ for its strict counterpart). This theorem has as an immediate consequence $AC^L$, and also the existence of an $L$-definable ordinal enumeration of $L$.
 
 \begin{center}
\vspace*{-4.3cm}
    \hspace{1cm}	\begin{tikzpicture}
        \def\nangle{60} \def\r{3} \def\mygreen{green!50!black}
            
            \path (0,0) -- (0,7) node[pos=-0.05]
            {$L_0=\varnothing$}
            ;
            
            \path (0,0) -- (0,7) 
            node[pos=0.22,sloped] {$\dots$}
            node[pos=0.42,sloped] {$\dots$}
            ;
            
        \draw[line width=1pt] (90+0.5*\nangle:0.2*\r) --  (90-0.5*\nangle:0.2*\r);
        \draw[line width=1pt] (90+0.5*\nangle:0.4*\r) --  (90-0.5*\nangle:0.4*\r);
        \draw[line width=1pt] (90+0.5*\nangle:0.75*\r) --  (90-0.5*\nangle:0.75*\r);
        \draw[line width=1pt] (90+0.5*\nangle:0.95*\r) --  (90-0.5*\nangle:0.95*\r);
        
        \draw (90-0.5*\nangle:0.2*\r) node [xshift=0.5cm] {$L_1$};
        \draw (90-0.5*\nangle:0.4*\r) node [xshift=0.5cm] {$L_2$};
        \draw (90-0.5*\nangle:0.75*\r) node [xshift=0.5cm] {$L_\omega$};
        \draw (90-0.5*\nangle:0.95*\r) node [xshift=0.7cm] {$L_{\omega+1}$};
        
        \draw[line width=2pt, cap=butt, rounded corners] (90+0.5*\nangle:1.1*\r) -- (0,0) node[yshift=1pt,fill=black,circle,minimum size=5pt, inner sep=0pt] {} -- (90-0.5*\nangle:1.1*\r);
        
        \draw (1.7, 3.1) node {$L$};
    \end{tikzpicture}
    \\
    \textit{L is stratified by its ordinal levels, which are \\ initial segments of the model (under $\leqslant_{L_\alpha}$).}
\end{center}
\vspace*{-0.5cm}
\begin{theorem}{\bf  (Condensation) (II.5.2 in  \cite{Dev})}
\label{bottLine}
If $\alpha$ is a limit ordinal and 
$X \prec_1 L_\alpha$, then there are unique $\pi$ and $\beta\leqslant\alpha$ such that
\begin{itemize}
    \item[\normalfont i)] $\pi: \langle X, \in\rangle \cong \langle L_\beta, \in\rangle$
    
    \item[\normalfont ii)] for transitive $Y \subseteq X, \, \pi\restriction Y = Id \restriction Y$
    
    \item[\normalfont iii)] $\pi(x) \leqslant_{L_\alpha} x$ for all $x \in X$
    
\end{itemize}
\end{theorem}

$A \prec_1 B$ means as usual that $A \subseteq B$ and both models satisfy the same $\Sigma_1$ formulas with parameters in $A$. Notice $\pi$ is just the Mostowski collapse, defined recursively as $\pi(y)=\{\pi(x) : x \in^X y\}$ (see Section I.9 in \cite{Kun}).

\vspace{0,2cm}

Condensation is arguably the most important tool for the study of constructibility. The following lemma, proved through Condensation, has as an immediate consequence $GCH^L$.
\vspace{-0,5cm}
\begin{lemma} {\bf (VI.4.6 in \cite{Kun2})}
For any infinite ordinal $\alpha$, $\mathcal{P}(L_\alpha)^L = \mathcal{P}(L_\alpha)\cap L \subseteq L_{\alpha^{+L}}$ \\ where $\alpha^{+L}$ denotes the smallest cardinal bigger than $\alpha$ according to L.
\end{lemma}

$(\alpha^{+n})^L$ denotes the ordinal operation $\alpha^{+L}$ iterated $n$ times, and $\mathcal{P}^n$ the powerset operation iterated $n$ times. We use throughout a fixed pairing function for the naturals, denoted by $J$, presupposed to be primitive recursive.

\section{Gaps of reals}

A main focus of our study will be the details of the generation process for the constructible reals. As presented in \cite{Gaps}, some results centered around recursion of Putnam in \cite{Put} lead to a deeper set theoretic study of the hierarchy levels where no new reals appear.
\begin{definition}
$\alpha$ is a \emph{gap ordinal} iff $(L_{\alpha+1} \setminus L_\alpha)\cap \mathcal{P}(\omega) = \varnothing$ \end{definition}

In the next results, by $F$ being an ordinal operation $\Sigma_1$ definable in $L$, we mean a function from $ON^n$ to $ON$ (for some $n$) such that there's a $\Sigma_1$ formula $\Phi$ satisfying $F(\alpha_1, \alpha_2, ..., \alpha_n) = \alpha$ iff $L \models \Phi (\alpha_1, \alpha_2, ..., \alpha_n, \alpha)$.

We first prove a necessary technical result. The original idea is found in a proof of L\'{e}vy (Theorem 36 in \cite{Lev}), but we instead employ simpler modern methods.
\vspace{-0,2cm}
\begin{lemma} {\bf (L\'{e}vy)}
Let F be an ordinal operation $\Sigma_1$ definable in L. If at least one of the $\alpha_i$ is infinite, then $$F(\alpha_1, \alpha_2, ..., \alpha_n) < \alpha_1^{+L} + \alpha_2^{+L} + ... + \alpha_n^{+L}$$
\begin{proof}
Since $L \models \exists\alpha(F(\alpha_1, \alpha_2, ..., \alpha_n) = \alpha)$, there's such an $\alpha \in L$. Take a limit $\beta$ such that $\alpha \in L_ \beta$, and also $\alpha_1, \ldots, \alpha_n \in L_\beta$. Say $\gamma = \alpha_1^{+L} + \alpha_2^{+L} + ... + \alpha_n^{+L}$. Inside $L$, apply the Downward L\"{o}wenheim-Skolem Theorem to find an elementary submodel $X$ of $L_\beta$ of cardinality smaller than $\gamma$, and containing $\alpha, \alpha_1, \ldots, \alpha_n$ (as well as all their elements). By Condensation $X$ is isomorphic to some $L_\xi$, with $\xi < \gamma$. And since the $\alpha_i$ will be fixed by the Mostowski collapse, $L_\xi \models \exists\alpha(F(\alpha_1, \alpha_2, ..., \alpha_n) = \alpha)$. But since $F(\alpha_1, \alpha_2, ..., \alpha_n) = \alpha$ is $\Sigma_1$ expressible, and hence upwards absolute, $\alpha$ is really the value of $F(\alpha_1, \alpha_2, ..., \alpha_n)$, and is below $\gamma$.
\end{proof}
\end{lemma}
\begin{center}
\vspace{0.2cm}
		\newcommand\irregularcircle[2]{% radius, irregularity
  \pgfextra {\pgfmathsetmacro\len{(#1)+rand*(#2)}}
  +(0:\len pt)
  \foreach \a in {10,20,...,350}{
    \pgfextra {\pgfmathsetmacro\len{(#1)+rand*(#2)}}
    -- +(\a:\len pt)
  } -- cycle
}
		\begin{tikzpicture}

		        \draw[rounded corners = 1mm, line width=2pt](0, 0) \irregularcircle {0.7cm}{1mm};
		        
		        \draw[line width = 2pt] (-4, -1.3) -- (-5.5, 1.3) -- (-2.5, 1.3) -- cycle;
		        
		        \draw[line width = 2pt, xshift=0.6cm] (3.5, -0.8) -- (4.5, 0.8) -- (2.5, 0.8) -- cycle;
		        
		        \draw (-4, -1.3) node[yshift=0pt,fill=black,circle,minimum size=4pt, inner sep=0pt]{};
		        
		        \draw (3.5, -0.8) node[yshift=0pt,fill=black,circle,minimum size=4pt, inner sep=0pt, , xshift=0.6cm]{};
		        
		        \draw (4.9, -0.5) node {$\bm{L_\xi}$};
		        \draw (1.1, -0.5) node {$\bm{X}$};
		        \draw (-2.9, -0.5) node {$\bm{L_\beta}$};
		        \draw (-1.95, +0.27) node {\begin{footnotesize}L\"owenheim\end{footnotesize}};
		        \draw (2, +0.27) node {\begin{footnotesize}Condensation\end{footnotesize}};
		        
		        \draw[->, line width=2pt] (-2.9, 0) -- (-0.9, 0);
		        \draw[->, line width=2pt] (1, 0) -- (3.2, 0);
		        
		\end{tikzpicture}\\ 
		\textit{In the previous and many following proofs, we apply the Downward \\ L\"owenheim-Skolem Theorem followed by the Condensation Theorem \\ to transfer some properties of a constructible level into a smaller one.}
		\vspace{0.2cm}
		\end{center}

Furthermore, $L_\gamma \models F(\alpha_1, \alpha_2, ..., \alpha_n) = \alpha$, since an identical argument shows the necessary witness for the $\Sigma_1$ formula will be $<\gamma$. Using this result, we can prove there are arbitrarily high big gaps in a strong sense.
\vspace{-0,3cm}
\begin{theorem} {\bf (1.3 in \cite{Gaps})}
Let F be an ordinal operation $\Sigma_1$ definable in L, and let $\beta_1$, \ldots, $\beta_m$ $\in \omega^L_1$. Then there are arbitrarily large $\alpha \in \omega^L_1$ such that
$$(L_{F(\alpha, \vv{\beta})} \setminus L_\alpha) \cap \mathcal{P}(\omega) = \varnothing$$
\begin{proof}
Take $\gamma \in \omega^L_1$. If $\alpha < \omega_2^L$, then by Lemma 1.8, $F(\alpha, \beta_1, \ldots, \beta_m) < \alpha^{+L} + \beta_1^{+L} + ... + \beta_m^{+L} \leqslant \omega^L_2$. Since all the constructible reals are in $L_{\omega^L_1}$ (Lemma 1.6), by the absoluteness of $F$,
$$L_{\omega^L_2}\models \exists\alpha(\alpha>\gamma \wedge (L_{F(\alpha, \vv{\beta})} \setminus L_\alpha) \cap \mathcal{P}(\omega) = \varnothing)$$
Inside $L$, we apply the L{\"o}wenheim-Skolem Downward Theorem to find an elementary countable submodel $X$ of $L_{\omega^L_2}$. By the absoluteness of the satisfaction predicate, $X$ really is an elementary submodel in $V$. By Condensation $X$ is isomorphic to some $L_\xi$, $\xi<\omega^L_1$ due to its cardinality in $L$. By the absoluteness of "being $L_\alpha$" for the hierarchy levels above $\alpha$ and since $L_\xi \equiv L_{\omega^L_2}$, we get the desired $\alpha$.
\end{proof}
\end{theorem}

\begin{center}
\vspace{-4cm}
    			\begin{tikzpicture}
        \def\nangle{60} \def\r{3} \def\mygreen{green!50!black}

            \path (0,0) -- (0,7) 
            node[pos=0.24,sloped] {$\dots$}
            ;
            
       \path (0,0) -- (0,7) 
            node[pos = 0.08, yshift=0pt,fill=black,circle,minimum size=3pt, inner sep=0pt] {}
            
            node[pos = 0.13, yshift=0pt,fill=black,circle,minimum size=3pt, inner sep=0pt] {}
            
            node[pos = 0.18, yshift=0pt,fill=black,circle,minimum size=3pt, inner sep=0pt] {}
           
            ;
       
        \draw[line width=1pt] (90+0.5*\nangle:0.75*\r) --  (90-0.5*\nangle:0.75*\r);
        
        \draw [->, line width= 2pt] (0,1.95)--(0,3);
        \draw (90-0.5*\nangle:0.75*\r) node [xshift=0.5cm, yshift=-0.15cm] {$L_{\omega_1^L}$};

        \draw[line width=2pt, cap=butt, rounded corners] (90+0.5*\nangle:1.1*\r) -- (0,0) node[yshift=1pt,fill=black,circle,minimum size=5pt, inner sep=0pt] {} -- (90-0.5*\nangle:1.1*\r);
        
        \draw (1.7, 3.1) node {$L$};
    \end{tikzpicture}\\
    \textit{As an immediate consequence (choosing F as the + 1 operation), there \\ are arbitrarily big gap ordinals below $\omega_1^L$. Also, since all constructible \\ reals are in $L_{\omega_1}^L$, all ordinals $\geqslant \omega_1^L$ will trivially be gap ordinals.}
    \vspace{0.2cm}
\end{center}

This can be generalized in certain directions stated in \cite{Gaps}. We present below the most general version and supply a proof.
\vspace{-0,2cm}
\begin{lemma}
$\mathcal{P}^n(\kappa)^L \subseteq L_{(\kappa^{+n})^L}$
\begin{proof}
$\mathcal{P}(\kappa)^L \subseteq \mathcal{P}(L_\kappa)^L \subseteq L_{\kappa^{+L}}$ (Lemma 1.6), and if $\mathcal{P}^n(\kappa)^L \subseteq L_{(\kappa^{+n})^L}$, then \\ $\mathcal{P}(\mathcal{P}^n(\kappa))^L \subseteq \mathcal{P}(L_{(\kappa^{+n})^L})^L \subseteq L_{(\kappa^{+{n+1}})^L}$ again.
\end{proof}
\end{lemma}

\begin{theorem} {\bf (Generalization of 1.9)}
Let F be an ordinal operation $\Sigma_1$ definable in L, s $\leqslant$ n positive integers, and $\beta_1$, \ldots, $\beta_m$ $\in (\kappa^{+s})^L$. Then there are arbitrarily large $\alpha \in (\kappa^{+s})^L$ such that
$$(L_{F(\alpha, \vv{\beta})} \setminus L_\alpha) \cap \mathcal{P}^n(\kappa) = \varnothing$$
\begin{proof}
Take $\gamma \in (\kappa^{+s})^L$. If $\alpha \in (\kappa^{+{n+1}})^L$, then $F(\alpha, \vv{\beta}) < \alpha^{+L} + \beta_1^{+L} + ... + \beta_m^{+L} \leqslant (\kappa^{+{n+1}})^L$. Since $\mathcal{P}^n(\kappa)^L \subseteq L_{(\kappa^{+n})^L}$ (previous lemma), and by the absoluteness of $F$,
$$L_{(\kappa^{+{n+1}})^L}\models \exists\alpha(\alpha>\gamma \wedge (L_{F(\alpha, \vv{\beta})} \setminus L_\alpha) \cap \mathcal{P}^n(\kappa) = \varnothing)$$
Inside $L$, we apply the L{\"o}wenheim-Skolem Downward Theorem to find an elementary submodel of $L_{(\kappa^{+{n+1}})^L}$ containing $\gamma$ and $\vv{\beta}$, and of cardinality smaller than $(\kappa^{+s})^L$. By the absoluteness of the satisfaction predicate, it really is an elementary submodel in $V$, and its cardinality is less than $(\kappa^{+s})^L$. By Condensation, it is isomorphic to some $L_\xi$, with $\xi<(\kappa^{+s})^L$ due to its cardinality. By the absoluteness of "being $L_\alpha$" for the hierarchy levels above $\alpha$ and since $L_\xi \equiv L_{(\kappa^{+{n+1}})^L}$, we get the desired $\alpha$.
\end{proof}
\end{theorem}

The result can also be trivially generalized by replacing the cardinal $\kappa$ by any constructible transitive set $a$, and consequently $(\kappa^{+s})^L$ by $(|a|^{+s})^L$.

Notice that $F$ is required to be definable in $L$, not just in $V$, otherwise $L$ might not prove $F$ is a function.\footnote{An extreme example of this can be obtained by choosing $F(\alpha)$ to be the least ordinal greater than $\alpha$ such that there exists a transitive model containing that ordinal and satisfying $ZFC$ and the existence of a mesurable cardinal. Then in a universe with an unbounded class of inaccessible cardinals plus a measurable cardinal this is a definable function, but in its $L$ none of its values will be defined.}

Notice also that the result is not straightforwardly generalizable beyond $\Sigma_1$ functions, since more complex functions might not be upwards absolute and thus Lemma 1.8 would fail. Consider for instance the function $F(\alpha) = |\alpha|^+$, which is $\Pi_1$ definable. Clearly there is no infinite $\alpha \in \omega_1^L$ such that $(L_{|\alpha|^+} \setminus L_\alpha) \cap \mathcal{P}(\omega) = \varnothing$, since $|\alpha|^+ = \omega_1$, and the reals of $L$ appear cofinally in the levels below $\omega_1^L$, which is $\leqslant \omega_1$.

\vspace{0,8 cm}

We now work towards proving Theorem 1.23, a characterization revealing the connection between the generation of reals and the model-theoretical structure of the hierarchy. For that we first need some technical results.

From now on we assume $\alpha$ infinite, since the finite case is trivial: no finite ordinal is a gap.

\begin{definition}
$\alpha$ \emph{starts a gap} iff $\alpha$ is a gap ordinal and $\forall\beta$$<$$\alpha((L_\alpha \setminus L_\beta)\cap \mathcal{P}(\omega) \neq \varnothing)$
\end{definition}

We first prove a useful result: if a new real appears, then it actually does so by a definition without parameters.

\vspace{-0,3cm}

\begin{lemma} {\bf (Lemma 1 in  \cite{Boo})}
If $\alpha$ is not a gap ordinal, then there is a real not in $L_\alpha$ which is definable without parameters over $L_\alpha$.
\begin{proof}
Suppose $\phi(x, a)$ defines a real over $L_\alpha$, with $a\in L_\alpha$ a parameter, and that this real is not in $L_\alpha$ (without loss of generality, we can suppose the definition requires only one parameter). Then by using the well-order $<_{L_\alpha}$ of $L_\alpha$ definable without parameters in $L_\alpha$, the formula
\begin{center}
$ \exists y (\phi(x, y) \wedge \forall w (\phi(w, y) \rightarrow w \in \omega) \wedge \neg \exists z \forall w (w \in z \leftrightarrow \phi(w, y)) \newline \wedge \forall y' (y' <_{L_\alpha}  y \wedge  \forall w (\phi(w, y) \rightarrow w \in \omega) \rightarrow \exists z \forall w (w \in z \leftrightarrow \phi(w, y'))))$
\end{center}
\noindent also defines a new real, and has no parameters. Indeed, it will be the new real defined by $\phi$ with the $<_{L_\alpha}$-least parameter (and we know at least one such parameter exists).
 \end{proof}\end{lemma}

The definability without parameters of any new real allows us to complete the following model theoretic argument.
\vspace{-0,2cm}
\begin{lemma} {\bf (8.1 in  \cite{Gaps})}
If $\alpha$ is not a gap ordinal, then $L_\alpha$ is pointwise definable.
\begin{proof}
By 1.2, $Def^{L_\alpha} \prec L_\alpha$. By Condensation, $Def^{L_\alpha} \cong L_\xi$ for some $\xi \leqslant \alpha$. But by the previous lemma there's a real definable without parameters over $L_\alpha$ not in $L_\alpha$, and since $\omega$ is definable without parameters in $L_\xi$ this real is also definable over $L_\xi$. So it belongs to $L_{\xi+1}$, and thus $\xi = \alpha$ and $Def^{L_\alpha} = L_\alpha$.
\end{proof}
\end{lemma}

We now introduce the concept of an arithmetical copy, originated in Boolos \cite{Boo}. The regular presence of arithmetical copies will have relevant consequences for the structure of the constructible hierarchy, as we see later.

\begin{definition}
An \emph{arithmetical copy} of L$_\alpha$ is a set E$_\alpha \subseteq \omega$ encoding through a fixed primitive recursive pairing J a subset of $\omega \times \omega$ isomorphic to $\in \: \restriction L_\alpha$. That is, letting $$Field(E_\alpha) = \{ n \in \omega \: | \: \exists m \in \omega ( J(n, m)\in E_\alpha \lor J(m, n)\in E_\alpha)\}$$ $$R(E_\alpha) = \{ \langle n, m \rangle \in \omega \times \omega \: | \: J(n, m) \in E_\alpha\}$$ there is an isomorphism $\langle L_\alpha, \in \rangle \cong \langle$Field$(E_\alpha), R(E_\alpha) \rangle$. \\ We sometimes just write $E_\alpha$ for $R(E_\alpha)$.
\end{definition}

If $L_\alpha$ is pointwise definable, an obvious argument enumerating the formulas defining its elements yields an arithmetical copy, as in the following lemma.
\vspace{-0,2cm}
\begin{lemma} {\bf (4.1 in  \cite{Gaps})}
If $L_\alpha$ is pointwise definable, then there is an arithmetical copy $E_\alpha$ of $L_\alpha$ belonging to $L_{\alpha+2}$, and the isomorphism witnessing that is also in $L_{\alpha+2}$.
\begin{proof}
Denote $$ \phi(n, a) \equiv L_\alpha \models \phi_n (a) \wedge \forall b (b \neq a \rightarrow \neg L_\alpha \models \phi_n (b))$$
$$ \Phi(n, a) \equiv \phi(n, a) \wedge \forall m < n (\neg\phi(m, a))$$
So $\Phi(n, a)$ states "$n$ is the least G\"{o}del number of a formula defining $a$". We can define
$$ E_\alpha = \{ J(n, m) \: | \: L_{\alpha+1} \models \exists a, b ( \Phi(n, a) \wedge \Phi(m, b) \wedge a\in b) \}$$
$$ \pi_\alpha = \{ \langle n, a \rangle \: | \: L_{\alpha+1} \models \Phi(n, a) \}$$
which thus belong to $L_{\alpha+2}$.
\end{proof}
\end{lemma}

But for Theorem 1.23 we need a stronger result: we need to see this copy actually belongs to $L_{\alpha+1}$. For this, a more fine-structural analysis by Boolos involving Skolem functions is needed, even if the general idea behind the proof remains the same. We require a flat pairing function.

\begin{definition}
A \emph{flat pairing function} is a definition of ordered pairs $P(x, y)$ (thus, preserving the property that for certain primitive recursive functions $\pi_1(P(x, y)) = x$ and $\pi_2(P(x, y)) = y$) that doesn't raise the constructibility rank. \\ That is, for infinite $\alpha$, $\forall x, y \in L_\alpha(P(x, y) \in L_\alpha)$.
\end{definition}
\vspace{-0,5cm}
\begin{lemma}
There is a flat pairing function.
\begin{proof}
Define $P(x, y) = x^0 \cup y^1$, where $x^0$ is obtained by replacing every natural number $n \in a \in x$ by its successor and adding 0 to $a$, and $y^1$ analogously but without adding 0.
\end{proof}
\end{lemma}

\begin{lemma} {\bf (Theorem 1 in  \cite{Boo})}
If $\alpha$ is not a gap ordinal, then there is an arithmetical copy $E_\alpha$ of $L_\alpha$ belonging to $L_{\alpha+1}$.
\begin{proof}
By Lemma 1.13, take a real $A \notin L_\alpha$ defined over $L_\alpha$ by $\phi(x)$ without parameters. Let $S$ be the Skolem hull closing $\omega$ in $L_\alpha$ under the Skolem functions for $\phi(x)$, $\neg\phi(x)$ and $V = L$ (that is, the Skolem functions for all of their existential subformulas). $S$ belongs to $L_{\alpha+1}$. Indeed, if $f_1, \ldots, f_m$ are these functions, each with $k_i$ variables, and $k = max_i\{k_i\}$, then

$ S = \{ x \in L_\alpha \: | \: L_\alpha \models \exists f, i (Function(f) \wedge Dom(f)=i+1 \wedge f(0)\in\omega \wedge \forall j<i$

$(f(j+1)\in\omega \lor \exists l_1\leqslant j \ldots l_k \leqslant j (f(j+1)=f_1(f(l_1), \ldots, f(l_{k_1})) \lor \ldots \lor$

$f(j+1) = f_m(f(l_1), \ldots, f(l_{k_m})))) \wedge f(i)=x   \}$

\noindent In case $\alpha$ is a successor ordinal, we encode all of these functions through iterated use of a flat pairing function, so that they all belong to $L_\alpha$ (previous lemma).

By construction, $S \prec_{V=L} L_\alpha$. This is actually the Tarski-Vaught proof of the L\"{o}wenheim-Skolem Theorem (see $\cite{Tar}$), and we can see it by induction on the complexity of the subformulas of $V=L$. For a $\Sigma_0$ formula it is immediate, since the relation $\in$ coincides in both $S$ and $L_\alpha$. For adding connectives which are not quantifiers, the induction step is trivial. And if $S\prec_{\phi} L_\alpha$, then $S\prec_{\exists x \phi} L_\alpha$, since if one such $x$ exists in $L_\alpha$ (for a certain choice of the parameters of $\phi$), then the Skolem function for this subformula has added it to $S$ by its construction.

Since $L_\alpha \models V=L$, we have $S \models V=L$, and thus $S \cong L_\xi$ for a certain $\xi$ (see for instance II.6.16 in \cite{Kun}), and the isomorphism can only be the Mostowski collapse because $L_\xi$ is transitive. $\xi \leqslant \alpha$, since otherwise $S \subseteq L_\alpha$ would contain an $\in$-chain of order type greater than $\alpha$ (which is impossible by induction on $\alpha$). And by the closure under $\phi(x)$ and $\neg \phi(x)$, again by the Tarski-Vaught proof of the L\"{o}wenheim-Skolem Theorem, $A$ is also definable over $S$, and thus $\xi = \alpha$.
As seen in the above definition, every member of $S$ either is a natural number or was added at least once by repeated applications of the Skolem functions. So by mirroring this construction using the natural numbers, we can define a recursive coding over $L_\alpha$ analogous to the above definition of $S$. Writing it out explicitly:

$ CODE(n, x) \equiv \exists f, g, i (Function(f) \wedge Function(g) \wedge Dom(f)=Dom(g)=i+1 \wedge f(0)\in\omega \wedge g(0)=J(0, f(0)) \wedge \forall j<i([f(j+1)\in\omega \wedge g(j+1)=J(0, f(j+1))] \lor \exists l_1\leqslant j \ldots l_k \leqslant j ([f(j+1)=f_1(f(l_1), \ldots, f(l_{k_1})) \wedge g(j+1) = J(1, f(l_1), \ldots, f(l_{k_1}))] \lor \ldots \lor [f(j+1) = f_m(f(l_1), \ldots, f(l_{k_m})) \wedge g(j+1) = J(m, f(l_1), \ldots, f(l_{k_m}))])) \wedge \\ f(i)=x \wedge g(i) = n)$

$Code(n, x) = CODE(n, x) \wedge \forall m < n (\neg CODE(m, x))$

\vspace{0,1cm}

By calling $\psi(x)$ the formula defining $S$ above, we can define 
$$ E_\alpha = \{ J(n, m) \: | \: L_\alpha \models \exists x, y (\psi(x) \wedge \psi(y) \wedge x\in y \wedge Code(n, x) \wedge Code(m, y))\}$$

and so $\langle L_\alpha, \in  \rangle \cong \langle S, \in \rangle \cong \langle Field(E_\alpha), E_\alpha \rangle$  and $E_\alpha \in L_{\alpha+1}$. \end{proof}
\end{lemma}

In fact, by using the flat pairing function, we can see the isomorphism $\pi_\alpha$ is also in $L_\alpha$, improving again on Lemma 1.16. The proof here was necessarily more intricate because we can't readily use a formula expressing $L_\alpha \models \phi_n$ without jumping up a level.

\begin{lemma} {\bf (Lemma 2.5 in  \cite{Gaps})}
If $\alpha$ starts a gap, then $\alpha$ is a limit ordinal.
\begin{proof}
Suppose $\beta+1$ starts a gap. We essentially define a new real by diagonalization.

$\beta$ is not a gap ordinal, so the previous proof applies to it. Defining $S$ and $Code$ as in that proof, define then the real\\ $Y = \{ n \in L_{\beta+1} \: | \: L_{\beta+1} \models \exists x, a_1, \ldots, a_k(x = \{ a \in S\cap \omega \: | \: \langle S, \in \rangle \models \phi_m (a_1, ..., a_k, a) \} \wedge \\ \indent \: \: \: \: a_1, \ldots a_k \in S \wedge Code(n_1, a_1) \wedge ... \wedge Code(n_k, a_k) \wedge n = J(m, n_1, \ldots, n_k) \notin x )\}$

This is well-defined because every $x$ with that form will belong to $L_{\beta+1}$ (that is, $L_{\beta+1}$ is right in building those sets). To see this, say $S = \{ b \in L_\beta \: | \: L_\beta \models \psi(b) \}$ as in the previous proof.
Now, if $x = \{ a \in S\cap \omega \: | \: \langle S, \in \rangle \models \phi_m (a_1, ..., a_k, a) \}$, then $x = \{ a \in L_{\beta} \: | \: L_\beta \models a \in \omega \wedge \phi_m' (a_1, ..., a_k, a) \}$, where $\phi_m'$ is the formula resulting from bounding the quantifiers in $\phi_m$ by the formula $\psi$.

Clearly $Y \in L_{\beta+2}$, and since $\beta+1$ is a gap ordinal, we have $Y \in L_{\beta+1}$. So for some $m$, $k$ and $b_1, ..., b_k \in L_\beta$, $Y = \{ x \in L_\beta \: | \: L_\beta \models \phi_m(b_1, ..., b_k, x) \}$. Then by applying the inverse of the Mostowski collapse $\pi$ (which fixes the natural numbers), $Y = \{ x \in S \: | \: \langle S, \in \rangle \models \phi_m (a_1, ..., a_k, x) \} $, where $\pi(a_i) = b_i$. And for certain $n_1, ..., n_k$, we have $Code(n_i, a_i)$. But then $J(m, n_1, \ldots, n_k) \in Y$ iff $J(m, n_1, \ldots, n_k) \notin Y$.
\end{proof}
\end{lemma}

The previous lemma grants us a shortcut to prove the following result (which can also be proved independently by similar reasoning to the proof of the previous lemma).

\vspace{-0,4cm}

\begin{lemma} {\bf (Lemma 2.4 in  \cite{Gaps})}
Let $\alpha$ be a gap ordinal. If $X \in \mathcal{P}(\omega) \cap L_\alpha$ and $X$ is a real well-ordering, then the type of $X$ is less than $\alpha$.
\begin{proof}
Suppose $\alpha$ is the least gap ordinal such that there's an $X \in \mathcal{P}(\omega)\cap L_\alpha$ coding a well-ordering of type $\geqslant \alpha$. If $\alpha$ weren't the start of a gap, then $X$ would already belong to the $\alpha' < \alpha$ starting the gap, contradicting leastness. So $\alpha$ starts the gap and is thus a limit by the previous lemma. But then, since $L_\alpha = \cup_{\beta<\alpha} L_\beta$, $X \in L_\beta$ for a certain $\beta < \alpha$, again contradicting leastness. \end{proof}
\end{lemma}

We now have the necessary tools to prove the central result.

\vspace{-0,3cm}

\begin{corollary}
If $\alpha$ starts a gap, then $L_\alpha \models V=HC$.
\begin{proof}
The non-gap ordinals $\beta$ are cofinal in $\alpha$ since it starts a gap. For any one of these non-gap ordinals, by Lemmas 1.14 and 1.16, $L_\alpha$ contains an isomorphism between $L_\beta$ and a subset of $\omega$.
\end{proof}
\end{corollary}

\vspace{-0,3cm}

\begin{theorem}
$\alpha$ starts a gap iff $L_\alpha \models ZFC-P + V=HC$
\begin{proof}
If $L_\alpha \models ZFC-P + V=HC$, then $\alpha$ is a gap ordinal. Otherwise, we'd have a set $\{ n\in L_\alpha \: | \: L_\alpha \models n \in \omega \wedge \phi(n) \} \in L_{\alpha+1}\setminus L_\alpha$, which would constitute a failure of Replacement in $L_\alpha$, since then $L_\alpha \models \forall x (x \neq \{ n\in\omega \: | \: \phi(n) \})$. And if $\alpha$ doesn't start the gap, then for some $\beta<\alpha$, $\mathcal{P}(\omega)^{L_\alpha} = \mathcal{P}(\omega)^{L_\beta} \in L_\alpha $, so $L_\alpha \models |\mathcal{P}(\omega)|>|\omega|$ by Cantor's theorem, contradicting $V=HC$.

If $\alpha$ starts a gap, by the previous lemma $L_\alpha \models V=HC$. All of the axioms except for Comprehension, Replacement and Choice are immediate from $\alpha$ being a limit.

For Comprehension, given $A\in L_\alpha$ and $\phi(x)$, choose a $\beta < \alpha$ with $A\in L_\beta$. By 1.14 and 1.16, the isomorphism $\pi_\beta$ between $L_\beta$ and a real belongs to $L_\alpha$. So if $\{ x \in A \: | \: L_\alpha \models \phi(x) \} \notin L_\alpha$, then we also have $\{ n \in \omega \: | \: L_\alpha \models \exists x\in A (\phi(x) \wedge \langle n, x\rangle \in \pi_\beta\} \notin L_\alpha$, contradicting $\alpha$ being a gap ordinal.

For Choice, because of the global order definable in $L_\alpha$, it is implied by Replacement. Indeed, given a family of disjoint nonempty sets, we will be able to construct by Replacement the set of all the $<_{L_\alpha}$-least elements of these sets, which will be a choice set.

To prove Replacement, take $A\in L_\alpha$ and $\phi(\bar{x}, y)$. We need to see there's a $\gamma < \alpha$ such that $\forall \bar{a} \in [A]^k (L_\alpha \models \phi(\bar{a}, y) \Rightarrow y\in L_\gamma)$.

Since $L_\alpha \models V=HC$, we choose an injection $f$ from $[A]^k$ to $\omega$ belonging to $L_\alpha$. Also, for every ordinal $\delta < \alpha$ choose the $<_{L_\alpha}$-least $W_\delta \subseteq \omega$ such that $\langle \delta, \in \rangle \cong \langle Field(W_\delta), W_\delta \rangle $ in $L_\alpha$ (at least one such real well-order exists by restricting the arithmetical copy $E_\delta \in L_{\delta+2}$ to the images of ordinals). We can thus define over $L_\alpha$ the following real, in which we code all well-orderings of order type $\delta+1$ for the $L_\delta$ in which a new one of the desired $y$ appears $$Z = \{ J(m, n) \: | \: L_\alpha \models \exists \bar{a} \in [A]^k \: \exists \delta, y (\phi(\bar{a}, y) \wedge y\in L_{\delta+1}\setminus L_\delta \wedge m=f(\bar{a}) \wedge n \in W_{\delta+1})\}$$

Since $\alpha$ is a gap ordinal, this real belongs to $L_\alpha$ and thus to some $L_{\beta+1}$ with $\beta + 1 < \alpha$. But then we can see that every $y$ such that $L_\alpha \models \exists \bar{a} \in [A]^k \phi(\bar{a}, y)$ will belong to $L_{\beta+1}$, and so $\beta+1$ will be our desired $\gamma$. Otherwise, by choosing such a $y_0 \in L_{\delta+1}\setminus L_\delta$ with $\delta > \beta$, we'd be able to construct from $Z$ a well-order of order type $\delta+1$ in $L_{\beta+1}$, contradicting Lemma 1.21. Indeed, if $L_\alpha \models \phi(\bar{a}_0, y_0)$, and $m_0 = f(\bar{a}_0)$, and $Z = \{ x \in L_\beta \: | \: L_\beta \models \psi(x)\}$, then $W_{\delta+1} = \{ n \in L_\beta \: | \: L_\beta \models \psi(J(m_0, n))\} \in L_{\beta + 1}$. \end{proof}
\end{theorem}

This last proof of Replacement amounts to the proof of a Reflection principle, as often Replacement and Reflection are closely linked. In fact, in \cite{Gaps} Marek and Srebrny detour through second-order arithmetic to prove the following principle, which they see implies $L_\alpha$ satisfying Replacement:

For any $\phi$, there are arbitrarily big $\beta\in\alpha$ such that, $\forall \bar{a} \in L_\beta$,
$$L_\beta \models \phi(\bar{a}) \:\:\: \mbox{iff} \:\:\:  L_\alpha \models \phi(\bar{a})$$

\vspace{0,4cm}
Notice in proving this central result we've made extensive use of the reals' capability to code information about the hierarchy or other sets. Indeed, given the extreme canonicity of $\omega$ (and most importantly that every natural is definable without parameters and present in all infinite levels of the hierarchy), the constructions in countable levels can eventually be replicated by a real. It is this versatility that ensures strong mirroring between the construction of reals and hierarchy levels, and why each one of these two processes can tell us much about the other.

\section{Gaps of finite order}

As already hinted at by Theorem 1.11, we can generalize the notion of gap to finite iterations of the power set operation. $J_n$ will now denote a fixed primitive recursive pairing function for $\mathcal{P}^n(\omega)$.

\begin{definition}
$\alpha$ is a \emph{k-gap ordinal} iff $(L_{\alpha+1} \setminus L_\alpha)\cap \mathcal{P}^k(\omega) = \varnothing$ \end{definition}

\begin{definition}
$\alpha$ \emph{starts a k-gap} iff $\alpha$ is a k-gap ordinal and\\ \phantom{Definition 1.29. a   aaaaaaaaaaaaaaaaa} $\forall \beta < \alpha ((L_\alpha\setminus L_\beta)\cap \mathcal{P}^k(\omega)\neq \varnothing)$
\end{definition}

\begin{definition}
An \emph{n-analytical copy} of L$_\alpha$ is a set E$_\alpha \subseteq \mathcal{P}^n(\omega)$ encoding through a fixed primitive recursive pairing $J_n$ a subset of $\mathcal{P}^n(\omega) \times \mathcal{P}^n(\omega)$ isomorphic to $\langle L_\alpha, \in \rangle$.
\end{definition}

Notice these definitions are slightly different (and not equivalent) to those of \cite{Gaps} in order to simplify notation.

Recall that $\mathcal{P}^0(\omega) = \omega$ and $\mathcal{P}^1(\omega) = \mathcal{P}(\omega)$. So the previous section corresponded to 1-gaps, and the arithmetical copies are the 0-analytical copies.

\vspace{0,4cm}

Since $\mathcal{P}^n(\omega) \subseteq \mathcal{P}^{n+1}(\omega)$, these definitions are equivalent to \\ $\forall n \leqslant k ((L_{\alpha+1} \setminus L_\alpha)\cap \mathcal{P}^n(\omega) = \varnothing)$ and $\forall\beta<\alpha\exists n \leqslant k ((L_\alpha \setminus L_\beta)\cap \mathcal{P}^n(\omega) \neq \varnothing)$ This exposes the inductive character that will facilitate generalizing the previous proofs. Since $\mathcal{P}^0(\omega) \subseteq \mathcal{P}^{k}(\omega)$, no finite ordinal is a $k$-gap, so we keep assuming $\alpha$ infinite.

\vspace{0,2cm}

Towards proving a stronger central theorem, we now generalize the previous results, providing only a sketch of the changes in the proofs where not obvious.

\vspace{-0,3cm}

\begin{lemma} {\bf (Generalization of 1.13)}
If $\alpha$ is not a k-gap ordinal, then there is an $x \in \mathcal{P}^{k}(\omega)$ not in $L_\alpha$ definable without parameters over $L_\alpha$.
\end{lemma}

\vspace{-0,3cm}

\begin{lemma} {\bf (Generalization of 1.14)}
For $k \geqslant 1$, if $\alpha$ is not a k-gap ordinal, then $L_\alpha$ is pointwise definable from $\mathcal{P}^{k-1}(\omega)^{L_\alpha}$.
\begin{proof}
We use $Def^{L_\alpha} (\mathcal{P}^{k-1}(\omega)^{L_\alpha})$ instead of just $Def^{L_\alpha}$ to ensure $L_\xi$ contains all of the members of $\mathcal{P}^{k-1}(\omega)^{L_\alpha}$ necessary for the definition to work and yield the same $x$.
\end{proof}
\end{lemma}

\vspace{-0,5cm}

\begin{lemma} {\bf{(Generalization of 1.16)}}
If $L_\alpha$ is pointwise definable from $\mathcal{P}^{k}(\omega)^{L_\alpha}$, then there is a k-analytical copy $E_\alpha$ of $L_\alpha$ belonging to $L_{\alpha+2}$, and the isomorphism witnessing that is also in $L_{\alpha+2}$.
\begin{proof}
We now have, for $a\in L_\alpha$ and $x\in \mathcal{P}^{k}(\omega)^{L_\alpha}$, $$ \phi(n, a, x) \equiv L_\alpha \models \phi_n (a, x) \wedge \forall b (b \neq a \rightarrow \neg L_\alpha \models \phi_n (b, x))$$
$$ \Phi(n, a, x) \equiv \phi(n, a, x) \wedge \forall m < n \forall y\in\mathcal{P}^{k}(\omega) (\neg\phi(m, a, y)) \wedge \forall y <_{L_\alpha} x (\neg\phi(n, a, y))$$
So $\Phi(n, a, x)$ states "$n$ is the least G\"{o}del number of a formula defining $a$ from a parameter $x$, and $x$ is the $<_{L_\alpha}$-least such parameter". Using $J_k$ a pairing function for $\mathcal{P}^k(\omega)$, and $i$ a fixed injection from $\omega$ into $\mathcal{P}^k(\omega)$, we define \\
$$ \hspace{-1,5cm} E_\alpha = \{ J_k(J_k(i(m), x), J_k(i(n), y)) \: | \: L_{\alpha+1} \models \newline $$ $$\: \: \: \: \: \: \: \: \: \: \: \: \: \: \: \: \: \: \: \: \: \: \: \: \: \: \: \: \exists x, y \in\mathcal{P}^k(\omega) \exists a, b ( \Phi(m, a, x) \wedge \Phi(n, b, y) \wedge a\in b \}$$
\indent \: \: \: \: \: \: \: \: \: \: \: $ \pi_\alpha = \{ \langle a, J_k(i(n), x) \rangle \: | \: L_{\alpha+1} \models \Phi(a, n, x) \}$
\end{proof}
\end{lemma}

\begin{lemma} {\bf (Generalization of 1.19)}
If $\alpha$ is not a k-gap ordinal and \\ $\mathcal{P}^{k-1}(\omega)^{L_\alpha} \in L_\alpha$, then there is a $(k-1)$-analytical copy $E_\alpha$ of $L_\alpha$ in $L_{\alpha+1}$.
\begin{proof}
By Lemma 1.27, there is an $A\in(L_{\alpha+1} \setminus L_\alpha)\cap \mathcal{P}^{k}(\omega)$ defined by $\phi(x)$ without parameters. Let $S$ be the hull closing $\mathcal{P}^{k-1}(\omega)^{L_\alpha}$ under the Skolem functions for $\phi(x)$, $\neg\phi(x)$ and $V=L$. As above, we need this to ensure $L_\xi$ contains all of the necessary members of $\mathcal{P}^{k-1}(\omega)$ to define $A$. The rest of the proof proceeds as before, using $<_{L_\alpha}$ for the members of $\mathcal{P}^{k-1}(\omega)$ instead of $<$ for the naturals, and $J_{k-1}$ and $i$ as in the previous lemma. \end{proof}
\end{lemma}

\begin{lemma} {\bf (Generalization of 1.20)}
If $\alpha$ starts a k-gap, then $\alpha$ is a limit ordinal.
\begin{proof}
Suppose $\beta+1$ starts a $k$-gap. If $x \in (L_{\beta+1} \setminus L_\beta)\cap \mathcal{P}^n(\omega)$, then clearly $\{ x \} \in (L_{\beta+2} \setminus L_{\beta+1})\cap \mathcal{P}^{n+1}(\omega)$, so by $\beta+1$ being a $k$-gap we can't have $n \leqslant k-1$. So $\beta$ is a $(k-1)$-gap, and thus $\mathcal{P}^{k-1}(\omega)^{L_{\beta+1}} = \mathcal{P}^{k-1}(\omega)^{L_{\beta}} \in L_{\beta+1}$. Then the former proof applies, using $J_{k-1}$ and $i$ and noticing the isomorphism between $S$ and $L_\alpha$ (which is the Mostowski collapse) fixes the members of $\mathcal{P}^{k-1}(\omega)$.\end{proof}
\end{lemma}

\begin{corollary}
If $\alpha$ starts a $k$-gap, then it isn't a $k'$-gap for any $k< k'$. \\ So in particular $\alpha$ can only start a gap of one finite order.
\begin{proof}
$\alpha$ is a limit and the non-$k$-gap ordinals are cofinal in $\alpha$. But then by defining over $L_\alpha$ the set of all members of $\mathcal{P}^{k}(\omega)$ in $L_\alpha$, we get a new member of $\mathcal{P}^{k+1}(\omega)$, and thus $\alpha$ isn't a $k'$-gap. 
\end{proof}
\end{corollary}

This new corollary provides insight on the generation of gaps of different orders. It implies that for a limit ordinal to be a $k$-gap, some limit below it already had to be part of its $(k-1)$-gap. So in particular, if $\alpha'$ starts a $k$-gap, then some $\alpha < \alpha'$ starts the 1-gap it's a part of, so between $\alpha$ and $\alpha'$ there are at least $k-1$ limits. Of course, by 1.11 we know $k$-gaps exist (and are cofinal in for instance $\omega_1^L$), so this is another way of interpreting the large length of some 1-gaps.

\begin{center}
\vspace*{-4cm}
    			\begin{tikzpicture}
        \def\nangle{60} \def\r{3} \def\mygreen{green!50!black}

            \path (0,0) -- (0,7) 
            node[pos=0.255,sloped] {$\dots$}
            ;
            \path (0,0) -- (0,7) 
            node[pos=0.355,sloped] {$\dots$}
            ;
            
       \path (0,0) -- (0,7) 
            node[pos = 0.2, yshift=0pt,fill=black,circle,minimum size=3pt, inner sep=0pt] {}
            node[pos=0.16] {1-gap}
            node[pos = 0.3, yshift=0pt,fill=black,circle,minimum size=3pt, inner sep=0pt] {}
            
            node[pos = 0.4, yshift=0pt,fill=black,circle,minimum size=3pt, inner sep=0pt] {}
           
            ;

        \draw[line width=2pt, cap=butt, rounded corners] (90+0.5*\nangle:1.1*\r) -- (0,0) node[yshift=1pt,fill=black,circle,minimum size=5pt, inner sep=0pt] {} -- (90-0.5*\nangle:1.1*\r);
        
        \path (0.3,0) -- (0.3,7)
            node[pos=0.4, xshift = 0.3cm]
            {3-gap} node[pos=0.3, xshift= 0.3cm]{2-gap};
    \end{tikzpicture}
    \\ \textit{For example, the only way for a 3-gap to appear will be for a 1-gap to extend \\ through two limit ordinals, the second of which will be the start of the 3-gap.}
    \begin{textblock*}{5cm}(11.8cm,5cm) % {block width} (coords) 
   Ordinals starting a 1-gap, \\ a 2-gap and a 3-gap
\end{textblock*}
\end{center}

\vspace{-0,3cm}

\begin{lemma} {\bf (Generalization of 1.21)}
Let $\alpha$ be a $k$-gap ordinal. If $X \in \mathcal{P}^{k}(\omega) \cap L_\alpha$ and $X$ codes a well-ordering through $J_{k}$, then the type of $X$ is less than $\alpha$.
\end{lemma}

\vspace{0,2cm}

Before proceeding to the generalization of the central result, let us note that the study of the generation of $\mathcal{P}^{k}(\omega)^L$ (that is, the study of $k$-gaps) is actually equivalent to the study of the generation of $\mathcal{P}(\aleph_{k-1})^L$. Indeed, the members of a higher cardinal can be used to code information in an equivalent way as we've been doing with $\omega$, only with a bigger cardinality of elements to choose from. And by the following Theorem 1.35, if $\alpha$ starts an $(n+1)$-gap, then $L_\alpha \models \mathcal{P}^n(\omega) \cong \aleph_{n}$. So this section can also be understood as a generalization of gaps to higher cardinals (although all of them below $\aleph_\omega^L$), as was also hinted at by 1.11, and the previous results can so be rephrased. Thus the following results come across as very natural.

\vspace{-0,3cm}

\begin{corollary} {\bf (Generalization of 1.22)}
\\ If $\alpha$ starts a $k$-gap, then $L_\alpha \models$ "there are k infinite cardinals"
\begin{proof}
A lower limit ordinal already was a $(k-1)$-gap, and thus \\ $\forall n \leqslant k-1 (\mathcal{P}^n(\omega)^{L_\alpha} \in L_\alpha)$, so $L_\alpha \models$ "$\aleph_{k-1}$ exists". As before, by Lemmas 1.28 and 1.29 every $L_\beta$ with $\beta<\alpha$ can be injected into $\aleph_{k-1}$.
\end{proof}
\end{corollary}

\vspace{-0,3cm}

\begin{theorem} {\bf (Generalization of 1.23)}
\\ $\alpha$ starts a $k$-gap iff $L_\alpha \models ZFC-P$ $+$ "there are k infinite cardinals"
\end{theorem}

The theory $ZFC-P$ $+$ \textit{"there are k infinite cardinals"} is in a sense as close to $ZFC$ as was possible: we want to add a subtheory of $ZFC$ to \textit{"there are k infinite cardinals"}, to ensure $L_\alpha$ is moderately right about how cardinal arithmetic works, but we can never have the Power Set axiom if there's a finite amount of cardinals (and thus a biggest cardinal). In fact, \textit{"there are k infinite cardinals"} is a sort of restricted Power Set axiom, stating that the first $k-1$ power sets of $\omega$ exist.

Of course, inside the theory $ZFC-P + V=L$ $+$ \textit{"there are k infinite cardinals"} we can define $(k+1)$-order constructible arithmetic. And in fact, by applying Theorem 2.1 of \cite{Zbi} to the constructible hierarchy, the levels modeling that theory are exactly those for which $L_\alpha \cap \mathcal{P}^{k}(\omega)$ is a model of this arithmetic. So we've actually seen how the generation of reals affects what these small models can say about the reals and successive constructions from them.

\section{Gaps of infinite order}

One might wonder whether there's an ordinal (for instance below $\omega_1^L$) that is a gap of every finite order. This might be understood intuitively as an $\omega$-gap, and so the generalization to infinite orders easily comes to mind. For this we first need to define what we mean by infinite iterations of the power set operation.

\begin{definition}
 
 $\mathcal{P}^0(\omega) = \omega$ 
 
 \phantom{aaaaaaa.aaaaa} $\mathcal{P}^{\beta+1}(\omega) = \mathcal{P}(\mathcal{P}^{\beta}(\omega))$

 \phantom{aaaaaaa.aaaaa} $\mathcal{P}^{\gamma}(\omega) = \bigcup_{\beta<\gamma} \mathcal{P}^{\beta}(\omega)$ for limit $\gamma$
\end{definition}

This is of course the only definition that makes the function $\mathcal{P}^{\alpha}(\omega)$ continuous on $\alpha$. The definitions of $\beta$-gap, start of a $\beta$-gap and $\beta$-analytical copy are as before. Notice that, for $\gamma < \gamma'$, $\mathcal{P}^\gamma(\omega) \subseteq \mathcal{P}^{\gamma'}(\omega)$, so the formulas in these definitions can again be rewritten as \\ \indent $\forall \gamma \leqslant \beta (L_{\alpha+1} \setminus L_\alpha)\cap \mathcal{P}^\gamma(\omega) = \varnothing$ and $\forall\delta<\alpha\exists \gamma \leqslant \beta ((L_\alpha \setminus L_\delta)\cap \mathcal{P}^\gamma(\omega) \neq \varnothing)$

\vspace{0,3cm}

By the reasoning after Lemma 1.33, the study of infinite order gaps is equivalent to the study of $\mathcal{P}(\kappa)$, where $\kappa$ can now be any cardinal (we work inside $L$ to simplify notation). Or also to the study of the finite order gaps over $\aleph_\gamma$ for $\gamma$ a limit. For instance, one might consider the gaps of order between $\omega$ and $\omega+\omega$ as the study of $\mathcal{P}^n(\aleph_\omega)$ for every $n \in \omega$.

In generalizing the previous results to these gaps, we will now have to consider also the case of $\beta$ being a limit ordinal.

\begin{lemma} {\bf (Generalization of 1.13)}
If $\alpha$ is not a $\beta$-gap ordinal, then there is an $x \in \mathcal{P}^{\beta}(\omega)$ not in $L_\alpha$ that is definable without parameters over $L_\alpha$.
\end{lemma}

\begin{lemma} {\bf (Generalization of 1.14)}
If $\alpha$ is not a $\beta$-gap ordinal, then $L_\alpha$ is pointwise definable from a certain $\mathcal{P}^{\gamma}(\omega)^{L_\alpha}$, for some $\gamma < \beta$.
\begin{proof}
For $\beta=\gamma+1$ we choose $\gamma$. For $\beta$ a limit, we choose a $\gamma < \beta$ such that $\alpha$ is not a $(\gamma+1)$-gap. Then the former proof applies.
\end{proof}
\end{lemma}

\begin{lemma} {\bf (Generalization of 1.16)}
If $L_\alpha$ is pointwise definable from $\mathcal{P}^{\gamma}(\omega)^{L_\alpha}$, then there is a $\gamma$-analytical copy $E_\alpha$ of $L_\alpha$ belonging to $L_{\alpha+2}$, and the isomorphism witnessing that is also in $L_{\alpha+2}$. \end{lemma}

\begin{lemma} {\bf (Generalization of 1.19)}
If $\alpha$ is not a $\beta$-gap ordinal and \\ $\forall \gamma < \beta( \mathcal{P}^\gamma(\omega)^{L_\alpha} \in L_\alpha)$, then there is a $\gamma$-analytical copy $E_\alpha$ of $L_\alpha$ in $L_{\alpha+1}$, for a certain $\gamma < \beta$.
\begin{proof}
As above, and by Lemma 1.37, there's a successor $\gamma + 1 \leqslant \beta$ with an $A\in(L_{\alpha+1} \setminus L_\alpha)\cap \mathcal{P}^{\gamma+1}(\omega)$ defined over $L_\alpha$ by $\phi(x)$ without parameters, so we argue as before with this $\gamma$. \end{proof}
\end{lemma}

\begin{lemma} {\bf (Generalization of 1.20)}
If $\alpha$ starts a $\beta$-gap, then $\alpha$ is a limit ordinal.
\begin{proof}
Suppose $\alpha+1$ starts a $\beta$-gap. For limit $\beta$, for some $\gamma < \beta$, $(L_{\alpha+1} \setminus L_\alpha)\cap \mathcal{P}^\gamma(\omega) \neq \varnothing$. But if $x \in (L_{\alpha+1} \setminus L_\alpha)\cap \mathcal{P}^\gamma(\omega)$, then clearly $\{ x \} \in (L_{\alpha+2} \setminus L_{\alpha+1})\cap \mathcal{P}^{\gamma+1}(\omega)$, contradicting $\alpha+1$ being a $\beta$-gap. For successor $\beta$, the former proof applies. \end{proof}
\end{lemma}

As before we can see that the start of a $\gamma$-gap is not a $(\gamma+1)$-gap. Notice that for starting a $\beta$-gap with $\beta$ limit, $\alpha$ has to be a limit of limits. Otherwise, the previous limit would not be a $\beta$-gap, so it wouldn't be a $\gamma$-gap for a certain $\gamma < \beta$, so $\alpha$ wouldn't be a $(\gamma+1)$-gap.

\begin{lemma} {\bf (Generalization of 1.21)}
Let $\alpha$ be a $\beta$-gap ordinal. If $X \in \mathcal{P}^{\gamma}(\omega) \cap L_\alpha$ for a $\gamma \leqslant \beta$ and $X$ codes a well-ordering through $J_{\gamma}$, then the type of $X$ is less than $\alpha$.
\end{lemma}

\begin{corollary} {\bf (Generalization of 1.22)}
For $\alpha \geqslant \beta$, if $\alpha$ starts a $\beta$-gap, then \\ $L_\alpha \models$ "the cardinals are the $\aleph_{\gamma}$ with $\gamma < \beta$"
\begin{proof}
As before, for $\gamma < \beta$ we can see the elements of $\mathcal{P}^\gamma(\omega)$ can't appear cofinally in the levels below $\alpha$. So $L_\alpha \models$ "$\aleph_{\gamma}$ exists". As before, by Lemmas 1.38 and 1.39 every $L_\beta$ can be injected into one of these $\aleph_{\gamma}$.
\end{proof}
\end{corollary}

Of course for $\beta = \gamma + 1$ we'll just have \textit{"$\aleph_\gamma$ is the biggest cardinal"}. But for $\beta$ a limit we won't have a biggest cardinal. This last situation is now compatible with the Power Set axiom but incompatible with Replacement, and so the central result changes slightly.

\begin{theorem} {\bf (Generalization of 1.23)}
For $\alpha \geqslant \beta$, $\alpha$ starts a $\beta$-gap iff either
\\ \indent i) $\beta$ is a successor and $L_\alpha \models ZFC-P$ $+$ "the cardinals are the $\aleph_{\gamma}$ with $\gamma < \beta$"
\\ \indent ii) $\beta$ is a limit and $L_\alpha \models ZFC-R +$ "the cardinals are the $\aleph_{\gamma}$ with $\gamma < \beta$"
\begin{proof}	
\textit{i)} is as before. For \textit{ii)}, if $L_\alpha \models ZFC-R +$ \textit{"the cardinals are the $\aleph_{\gamma}$ with $\gamma < \beta$"}, then $\alpha$ is a $\beta$-gap. Otherwise, suppose $\phi(x)$ defines without parameters over $L_\alpha$ a new member of $\mathcal{P}^{\gamma+1}(\omega)$, for a $\gamma < \beta$. This constitutes a failure of Comprehension in $L_\alpha$, since it is the set $\{ x \in \mathcal{P}^{\gamma}(\omega)^{L_\alpha} \: | \: L_\alpha \models \phi(x) \}$, and $\mathcal{P}^{\gamma}(\omega)^{L_\alpha} \in L_\alpha$.

And if $\alpha$ doesn't start the gap, then for some $\delta < \alpha$, $\forall \gamma < \beta (\mathcal{P}^\gamma(\omega)^{L_\alpha} = \mathcal{P}^\gamma(\omega)^{L_\delta})$ and so $\cup_{\gamma < \beta} \mathcal{P}^\gamma(\omega)^{L_\alpha} \in L_{\delta+1}$ contradicting the non-existence of $\aleph_\beta$ in $L_\alpha$.

For Choice, given a family of disjoint non-empty sets in $L_\alpha$, it belongs to some $L_\delta$ with $\delta<\alpha$. But there'll be a certain $\delta \leqslant \delta' < \alpha$ such that $\delta'$ starts a $\gamma'$-gap (for some successor $\gamma' < \beta$), since the non-$\beta$-gap ordinals are cofinal in $\alpha$ (because it starts the gap) and also every $\mathcal{P}^\gamma(\omega)^{L_\alpha}$ belongs to a certain level below $\beta$. By \textit{i)}, it will satisfy $L_{\delta'} \models ZFC-P$, and so a choice set for the family will belong to it.

For Power Set, the existence of the power of the $\aleph_\gamma$ ensures the existence of the power of every set. Indeed, suppose the subsets of a certain $x \in L_\alpha$ were cofinal in $L_\alpha$. We know $L_\alpha$ proves $x$ isomorphic to (a subset of) a $\gamma$-analytical copy $E_{\delta}$, for a certain $\gamma < \beta$. Thus, since thanks to the isomorphism they are mutually definable, the subsets of $E_{\delta}$ would also be cofinal in $L_\alpha$, but these are members of $\mathcal{P}^{\gamma+1}(\omega)^{L_\alpha}$, so that would contradict $\alpha$ being a $\beta$-gap. So the subsets of $x$ are not cofinal in $L_\alpha$ and thus $\mathcal{P}(x)^{L_\alpha}$ will be constructed at a certain level. \end{proof}
\end{theorem}

Going back to the question opening this section, the answer is positive by further generalizing Theorem 1.9, and we see our formulation of $\beta$-gaps is natural to express results such as this one.

\begin{theorem} {\bf (Generalization of 1.9)}
For any $\beta$ and $\gamma \geqslant 1$ such that $\beta < \omega_\gamma^L$, there are arbitrarily big $\beta$-gap ordinals below $\omega_\gamma^L$.
\begin{proof}
For $\gamma > \beta $ it is immediate, since $\mathcal{P}^\beta(\omega)^L \subseteq L_{\omega_\beta^L}$. For $\gamma \leqslant \beta$, take $\delta \in \omega_\gamma^L \subseteq \omega_\beta^L$. Then again by $\mathcal{P}^\beta(\omega)^L \subseteq L_{\omega_\beta^L}$
$$L_{\omega_{\beta+1}^L}\models \exists\alpha(\alpha>\delta \wedge (L_{\alpha+1} \setminus L_\alpha) \cap \mathcal{P}^\beta(\omega) = \varnothing)$$
Inside $L$, we apply the L{\"o}wenheim-Skolem Downward Theorem to find an elementary submodel of $L_{\omega_{\beta+1}^L}$ containing $\delta$ and $\beta$, and of cardinality smaller than $\omega_\gamma^L$ (possible because $\gamma \geqslant 1$). By Condensation, it is isomorphic to some $L_\xi$, with $\xi<\omega_\gamma^L$, and we get the desired $\alpha$.
\end{proof}
\end{theorem}

The only change in the proof is including $\beta$ as a parameter to ensure $\mathcal{P}^\beta(\omega)$ is correctly defined, that is, $\beta$ collapses to itself, and that forces the requirement $\beta < \omega_\gamma^L$. This premise is necessary, since for instance there isn't an $\omega_1$-gap ordinal below $\omega_1^L$. Indeed, if $\alpha_{\omega_1}$ started such a gap, by the reasoning after Lemma 1.41 there'd be an ordinal $\alpha_\gamma$ starting the $\gamma$-gap to which $\alpha_{\omega_1}$ belongs, for every $\gamma < \omega_1$. But $\gamma < \gamma'$ implies $\alpha_\gamma < \alpha_{\gamma'}$, and thus there would be an uncountable (according to $L$) number of ordinals below $\alpha_{\omega_1} < \omega_1^L$.

Just as in Theorem 1.9, the previous theorem can be generalized to an ordinal operation $\Sigma_1$ definable in $L$ (replacing the $+1$ operation) to obtain longer gaps.

\section{Lengths of gaps}

Returning to the results of \cite{Gaps}, we study now more closely the exact lengths of these gaps, and find striking regularity. We present the results generalized to $\beta$-gaps straight away.

\begin{definition}
$\alpha$ \emph{starts a $\beta$-gap of length $\rho$} iff it starts a $\beta$-gap,  \\ $(L_{\alpha + \rho} \setminus L_\alpha)\cap \mathcal{P}^\beta(\omega) = \varnothing$ and $(L_{\alpha + \rho + 1} \setminus L_\alpha)\cap \mathcal{P}^\beta(\omega) \neq \varnothing$
\end{definition}

Let us enumerate the beginnings of $\beta$-gaps as $\alpha_0$, $\alpha_1$, ..., $\alpha_\xi$, ... \\ Provided $\beta < \omega_\beta^L$, there are $\omega_\beta^L$ such beginnings, because there are arbitrarily big $\beta$-gap and non-$\beta$-gap ordinals in $\omega_\beta^L$.

When on the contrary $\beta=\omega_\beta^L$ (for instance, when $\beta = \omega_{\omega_{\omega_{._{._{.}}}}}^L$), by the reasoning ending the previous section there are no $\beta$-gaps below $\beta$. But of course all ordinals $\geqslant \omega_\beta^L$ will be $\beta$-gaps. So these surprisingly present a univocal cut-off point.

\begin{theorem} {\bf (Generalization of 4.2 in \cite{Gaps})}
The ($\xi$+1)st $\beta$-gap is of length 1, \\ provided $\beta \leqslant \alpha_\xi < \omega_\beta^L$.
\begin{proof}
$L_{\alpha_{\xi+1}}$ is pointwise definable from $L_{\alpha_\xi}\cup\{L_{\alpha_\xi}\}$. Indeed, $Def^{L_{\alpha_{\xi+1}}}(L_{\alpha_\xi}\cup\{L_{\alpha_\xi}\}) \prec L_{\alpha_{\xi+1}}$ by a Skolem hull argument, and so it is isomorphic to some $L_\gamma$ by Condensation, with $\gamma > \alpha_\xi$. But $\gamma$ must start a $\beta$-gap because it satisfies the theory corresponding to the starts of $\beta$-gaps by Theorem 1.44, so $\gamma = \alpha_{\xi+1}$.
In $L_{\alpha_{\xi+1}}$, the set $L_{\alpha_\xi}\cup\{L_{\alpha_\xi}\}$ is injectible into $\mathcal{P}^\gamma(\omega)$ for a certain $\gamma < \beta$ (that is, of cardinality $\aleph_\gamma$ or less), so $L_{\alpha_{\xi+1}}$ is also pointwise definable from $\mathcal{P}^\gamma(\omega)$, and by Lemma 1.39 we are done.
\end{proof}
\end{theorem}

Thus, again if $\beta < \omega_\beta^L$, there'll be $\omega_\beta^L$ $\beta$-gaps of length 1. But the result generalizes for any length in the following theorem, which is our main tool for dealing with lengths.

\vspace{-0.3cm}

\begin{theorem} {\bf (Generalization of 4.4 in \cite{Gaps})}
The first $\beta$-gap of length $\geqslant \rho$ starting (strictly) higher than a given $\gamma \geqslant \rho$ exists and is of length $\rho$, provided $\beta \leqslant \gamma < \omega_\beta^L$.
\begin{proof}
Say $\alpha$ starts said gap, and suppose $\rho = \delta+1$. We will build a new real in $L_{\alpha+\delta+2}$. Consider $M := Def^{L_{\alpha+\delta}}(\gamma\cup\{\gamma\})$. By a Skolem hull argument $M \prec L_{\alpha+\delta}$. By Condensation it is isomorphic to some $L_\xi$ through the Mostowski collapse, with $\xi \leqslant \alpha+\delta$. Clearly $\gamma\cup\{\gamma\}\subseteq L_\xi$, and also $\alpha \in M$, since it is definable as the $<$-least ordinal starting a $\beta$-gap of length $\geqslant \rho$ above $\gamma$. Let $\bar{\alpha}$ be the collapse of $\alpha$. It's greater than $\gamma$ (since the ordinals $\leqslant \gamma$ collapse to themselves), and starts a $\beta$-gap of length $\geqslant \rho$, so $\bar{\alpha} \geqslant \alpha$ by definition of $\alpha$, and thus $\bar{\alpha} = \alpha$. Since $L_{\alpha+\delta} \models \forall \delta'<\delta(\alpha+\delta'$ exists$)$, we have $\xi = \alpha+\delta$. Thus $L_{\alpha+\delta}$ is pointwise definable from $\gamma \cup \{\gamma\}$, which is injectible into $\mathcal{P}^{\gamma'}(\omega)$ for a certain $\gamma' < \beta$ because of the theory $L_\alpha$ satisfies, so again by Lemma 1.39 we are done.

Now suppose that $\rho$ is a limit. We construct a $\beta$-analytical copy $E_{\alpha+\rho} \in L_{\alpha+\rho+1}$. We choose a finite number of sentences that guarantee every well-founded model of them (containing $\gamma$) to be isomorphic to $L_{\alpha+\rho}$. These sentences are extensionality, $V=L$ and the sentence
\\ \indent $\exists \mu [\forall x(x \in \mathcal{P}^\beta(\omega) \rightarrow x \in L_\mu) \\ \indent \wedge \forall \nu \in \mu \exists x \in \mathcal{P}^\beta(\omega) (x \in L_\mu \setminus L_\nu) \\ \indent \wedge \forall \nu \in \rho (\mu + \nu$ exists$)]$

\vspace{0,2cm}

As in Lemma 1.19, we take the Skolem hull of $\gamma \cup \{ \gamma \}$ under the Skolem functions of these sentences in $L_{\alpha+\rho}$, which is isomorphic to $L_{\alpha+\rho}$, belongs to $L_{\alpha+\rho+1}$, and whose construction can be coded by members of $\mathcal{P}^{\gamma'}(\omega)$ for a certain $\gamma' < \beta$ (into which $\gamma \cup \{ \gamma \}$ is injectible). \end{proof}
\end{theorem}

Notice $\gamma < \omega_\beta^L$ is trivially required since otherwise there are no starts of $\beta$-gaps above it. And we need $\beta \leqslant \gamma$ to ensure the definition of $\alpha$ in the successor case and the last sentence in the limit case can be formulated.

\vspace{-0,3cm}

\begin{corollary}
If $\beta < \omega_\beta^L$, then for every $\rho \in \omega_\beta^L$ there are $\omega_\beta^L$ gaps of length $\rho$.
\end{corollary}

\vspace{-0,3cm}

\begin{lemma} {\bf (Generalization of 4.7 in \cite{Gaps})}
If $\alpha$ starts a $\beta$-gap of length $> 1$, the beginnings of $\beta$-gaps of length 1 are cofinal in $\alpha$, provided one of these beginnings is $\geqslant \beta$.
\begin{proof}
Otherwise, consider their supremum $S < \alpha$. Suppose there's an $S \leqslant \alpha' < \alpha$ beginning a $\beta$-gap of length $> 1$. Then by 1.47, after $\alpha'$ and before $\alpha$ there's a $\beta$-gap of length 1, contradicting $S$ being the supremum. So there are no $\beta$-gaps above (or equal to) $S$ and below $\alpha$.

Then by an argument as that of Theorem 1.47, $L_\alpha$ is pointwise definable from $L_S\cup\{L_S\}$ (because it's the first $\beta$-gap level higher than $S$), which is injectible into a certain $\mathcal{P}^{\gamma'}(\omega)$, so as before by Lemma 1.39 $\alpha$ starts a $\beta$-gap of length 1.
\end{proof}
\end{lemma}

\begin{theorem} {\bf (Generalization of 4.8 in \cite{Gaps})} If $\alpha$ starts a $\beta$-gap of length $\rho < \alpha$, \\ then for each $\sigma < \rho$, $sup\{\delta < \alpha \: | \: \delta$ starts a $\beta$-gap of length $\sigma \} = \alpha$, \\ provided any one of these beginnings is $\geqslant \beta$.
\begin{proof}
By the above lemma this is true for $\sigma = 1$. Consider now only these beginnings of $\beta$-gaps of length 1 which are above $\rho$. For each $\sigma < \rho$, by Theorem 1.48 the first $\beta$-gap of length $\geqslant \sigma$ starting higher than any one of the given $\beta$-gaps of length 1 is of length $\sigma$. So between every $\beta$-gap of length 1 and $\alpha$ there's a $\beta$-gap of length $\sigma$.
\end{proof}
\end{theorem}

This result exposes extreme regularity in the lengths of gaps, and the slow hierarchical building of them: a gap of a given length can only appear as the limit of many other gaps of smaller length.

\begin{center}
        \begin{tikzpicture}
        \def\nangle{70} \def\r{7.2} \def\mygreen{green!50!black}

            \path (0,0) -- (0,7)
            node[pos = 0.12, yshift=0pt,fill=black,circle,minimum size=3pt, inner sep=0pt] {}
            node[pos = 0.2, yshift=0pt,fill=black,circle,minimum size=3pt, inner sep=0pt] {} 
            node[pos = 0.28, yshift=0pt,fill=black,circle,minimum size=3pt, inner sep=0pt] {} 
            node[pos=0.16,sloped] {\tiny $\dots$}
            node[pos = 0.4, yshift=0pt,fill=black,circle,minimum size=4pt, inner sep=0pt] {} 
            node[pos = 0.48, yshift=0pt,fill=black,circle,minimum size=3pt, inner sep=0pt] {} 
            node[pos = 0.56, yshift=0pt,fill=black,circle,minimum size=3pt, inner sep=0pt] {} 
            node[pos = 0.64, yshift=0pt,fill=black,circle,minimum size=3pt, inner sep=0pt] {} 
            node[pos = 0.76, yshift=0pt,fill=black,circle,minimum size=4pt, inner sep=0pt] {}
            node[pos = 0.96, yshift=0pt,fill=black,circle,minimum size=5pt, inner sep=0pt] {}
            
            node[pos=0.24,sloped] {\tiny $\dots$}
            node[pos=0.35,sloped] {$\dots$}
            node[pos=0.7,sloped] {$\dots$}
            node[pos=0.44,sloped] {\tiny $\dots$}
            node[pos=0.52,sloped] {\tiny $\dots$}
            node[pos=0.6,sloped] {\tiny $\dots$}
            
            node[pos=0.87,sloped] {\Huge $\dots$}
            
            ;
            
            \path (0.3,0) -- (0.3,7) node[pos=0.28]
            {\small 1} node[pos=0.12, xshift = 0cm] {\small 1} node[pos=0.2, xshift = 0cm] {\small 1} node[pos=0.4, xshift = 0cm]
            {2} node[pos=0.56, xshift= 0cm]
            {\small 1} node[pos=0.48, xshift= 0cm]
            {\small 1} node[pos=0.64, xshift= 0cm]
            {\small 1}  node[pos=0.76, xshift = 0cm]
            {2} node[pos=0.96, xshift= 0cm]
            {\large 3};

        \draw[line width=2pt, cap=butt, rounded corners] (90+0.5*\nangle:1.1*\r) -- (0,0) node[yshift=1pt,fill=black,circle,minimum size=5pt, inner sep=0pt] {} -- (90-0.5*\nangle:1.1*\r);
    \end{tikzpicture} \\
        \begin{textblock*}{5cm}(10.1cm,10.7cm) % {block width} (coords) 
   Lengths of \\ the gaps started 
\end{textblock*}
    \vspace{0.2cm}
    \textit{For example, a gap of length 3 is the limit of gaps of \\ length 2, and thus the limit of limits of gaps of length 1.}
\end{center}

\section{An application}

Consider the following Friedman-Tomasik theorem of [23].

\begin{theorem} {\bf (II.7.3A-3E in \cite{Dev})}
There are $\omega_1^L$ theories of sets $L_\alpha$.
\end{theorem}

That is, if $\Sigma = \{ \{ \phi \: | \: L_\alpha \models \phi \} \} _{\alpha \in ON}$, then $L \models |\Sigma| = \aleph_1$. Even if $V \neq L$, this will of course yield an actual bijection between $\Sigma$ and $\omega_1^L$.

This result is usually proved by a diagonalization argument, but thanks to the link between gap ordinals and pointwise definability explored earlier, phrasing it in terms of gaps will make the generalization easier.

\begin{proof} Since $\Sigma$ can be injected into $\mathcal{P}(\omega)$ inside $L$ (by using G\"{o}del numbers), its cardinality in $L$ is $\leqslant \aleph_1$.
On the other hand, there are $\omega_1^L$ non-1-gap ordinals (they are cofinal in $\omega_1^L$). By Lemma 1.14, their corresponding levels are pointwise definable, and thus have no proper elementary submodels. Indeed, if $M \prec L_\alpha$, and $\phi$ defines $x\in L_\alpha$ without parameters over $L_\alpha$, then it also does for $M$. Otherwise it would define a different $y\in M$, and thus $M$ and $L_\alpha$ would disagree about $\phi(y)$.

But then any two of these levels have different theories, since otherwise they'd also be isomorphic (by identifying elements with the same definition). And that's not possible since then one would be a proper elementary submodel of the other. So we've also seen the cardinality of $\Sigma$ in $L$ is $\geqslant \aleph_1$. \end{proof}

\begin{definition}
We call the theory of a model of the form $$ \langle L_\alpha, \in, a \rangle _{a \in \mathcal{P}^\beta(\omega)\cap L_\alpha}$$
a $(\beta+1)$st-order analysis theory.
\end{definition}

\noindent The $a$ are of course used as constants, and thus the language has cardinality $\aleph_\beta^L$.

\begin{theorem} {\bf (Generalization of 1.52)}
\\ There are $\omega_{\beta+1}^L$ $(\beta+1)$st-order analysis theories of sets $L_\alpha$.
\begin{proof}
There are $\leqslant \omega_{\beta+1}^L$ such theories, again by using G\"{o}del numbers and the cardinality of the language. On the other hand, there are $\omega_{\beta+1}^L$ non-$(\beta+1)$-gap ordinals above $\omega_\beta^L$. For every such $\alpha \geqslant \omega_\beta^L$, by Lemma 1.38, $L_\alpha$ is pointwise definable from $\mathcal{P} ^\beta (\omega) ^{L_\alpha} = \mathcal{P} ^\beta (\omega) ^L$, and thus have no proper elementary submodels containing $\mathcal{P} ^\beta (\omega) ^L$ (by an argument as above). So they are pairwise non-isomorphic as before, and have different $(\beta+1)$st-order analysis theories. \end{proof}
\end{theorem}

Given the tight link between gaps and fundamental model-theoretic concepts like Condensation and pointwise definability, presumably many other results in constructibility can be meaningfully rephrased in terms of gaps. In that direction, Marek and Srebrny observe that appropriate generalizations of the notion of a gap will correspond to a diverse range of set theories, including statements about the existence of certain (small) ordinals or cardinals.

\newpage

\chapter{Pathologies in L}

Unlike other objects of Set Theory, the reals aren't well-behaved in $L$. We'll see some of the model's properties have pathological consequences for the structure of the reals. These might seem reasonable arguments against the adoption of the Axiom of Constructibility $V=L$, or even against $CH$.

\section{A well-ordering of the reals}

It proves hard to imagine how a well-ordering of the reals (or of any set of uncountable cardinality) might look like. This is nothing but another instance of our intuitions about continuous and discrete objects colliding. The continuity of the real line, that is, its dense order and closure under converging infinite sequences, is the very reason why it seems such a natural and necessary object in the first place. So it is not surprising that the idea of an order on it with radically different properties should seem alien.

Then of course, it does seem plausible that any set of discrete objects will be easily well-orderable: we just choose an object at each step for a transfinite amount of times. And (without delving into serious philosophical dispute) transfinite iteration is usually regarded as a more intuitively plausible generalization, a natural extension of the obvious induction and recursion principles. Since thanks to the foundational power of Set Theory all mathematical objects are sets, they will all be plausibly well-orderable. If this apparent plausibility doesn't transfer that easily to the reals it is precisely because we have a hard time picturing them as a discrete set (that is, because of the tension between their intensional nature and extensional underpinning mentioned earlier).

The previous paragraph gives just a sketch of why the Axiom of Choice implies the Well-ordering theorem (every set is well-orderable). The apparent plausibility mismatch between these two principles is captured in Jerry Bona's famous quote \cite{Kra}: \\
\indent \textit{"The Axiom of Choice is obviously true, \\ \indent the Well-ordering theorem is obviously false; \\ \indent and who can tell about Zorn's Lemma?"}

\vspace{0,3cm}

The Axiom of Choice is not widely accepted only because of its plausibility. Within a set theoretic framework, the Axiom of Choice is required for many fundamental results and constructions, especially regarding other mathematical branches like topology and geometry, more closely related to our intuitive understanding of space and the continuum.

So Choice is ever present, and thus in any of the many conceivable different universes in which it holds a well-ordering of the reals will actually exist, our ability to picture it notwithstanding. Whether this well-ordering is definable by a formula, though, and the complexity of this definition, might vary across models of $ZFC$ (that is, across extensions of the theory $ZFC$), as we'll see in the next chapter.

But for now, it is enough to notice that inside $L$ (that is, under the assumption $V=L$) there is such a formula, thanks to its definable global well-order. And in fact, keeping track of its detailed structure as in the reference for Theorem 1.4, we see that the well-ordering of the reals induced by the restriction of the global well-ordering of $L$ has complexity both $\Sigma_2^1$ and $\Pi_2^1$, and thus $\Delta^1_2$, in the analytical hierarchy. That is, with at most two alternate quantifiers, which quantify only over the reals (thus yielding a sentence of second-order arithmetic).

\section{Non-measurability and non-Baireness}

An argument (probably originally due to Sierpi\'{n}ski) centered around Fubini's theorem from mathematical analysis shows that any well-order of the reals is a subset of $\mathbb{R}^2$ that is not Lebesgue measurable.

\begin{definition}
A \emph{null subset of $\mathbb{R}^n$} is one with Lebesgue measure 0. \\ A \emph{co-null set} is one with null complement.
\end{definition}

\vspace{-0.3cm}

\begin{theorem} {\bf (Fubini's Theorem for null sets)}
Suppose $W \subseteq \mathbb{R}^2$ is measurable, and for $x \in \mathbb{R}$ let  $W_x = \{ y\: | \: (x, y) \in W\}$. Then $W$ is null iff $\{ x \in \mathbb{R} \: | \: W_x$ is null in $ \mathbb{R}\}$ is co-null in $\mathbb{R}$.
\end{theorem}

\vspace{-0,3cm}

\begin{theorem} {\bf (1.2 in \cite{Cai})}
Any well-ordering of a non-null set of reals is not Lebesgue measurable.
\begin{proof}
 Towards a contradiction, let $\lambda$ be the least ordinal such that there is a non-null $S \subseteq \mathbb{R}$, and an enumeration $\langle r_\alpha \: | \: \alpha < \lambda \rangle$ of $S$ such that $W = \{ (r_\alpha, r_\beta) \: | \: \alpha < \beta < \lambda \} \subseteq \mathbb{R}^2$ is measurable. Let $S_\alpha = \{ r_\beta \: | \: \beta < \alpha \}$ and $S^\alpha = \{ r_\beta \: | \: \alpha < \beta \}$. Also, for $x \in \mathbb{R}$, $r^{-1}(x)$ is the $\alpha$ such that $x = r_\alpha$.
 
 Let's see $S$ is measurable. Indeed, for almost all $x \in S$ and almost all $y \in S$, both $S^{r^{-1}(x)}$ and $S_{r^{-1}(y)}$ are measurable, because otherwise $W$ wouldn't be measurable. And since $S = S_{\gamma+1}\cup S^\alpha$ for any $\alpha \leqslant \gamma < \lambda$, $S$ must be measurable.
 
 Now let's find a $\gamma < \lambda$ such that $S_\gamma$ is non-null and measurable, and we will be done by contradicting the minimality of $\lambda$.
 
 Since for almost all $y \in S$ $S_{r^{-1}(y)}$ is measurable, if there is no such $S_\gamma$ then almost all of them are null, and thus the measure of $W$ is 0. But on the contrary, for almost all $x \in \mathbb{R}$, $S^{r^{-1}(x)}= S \setminus (S_{r^{-1}(x)} \cup \{x\})$ has positive measure (since $S$ is not null, but almost all $S_{r^{-1}(x)}$ are), so $W$ can't have measure 0.
\end{proof}
\end{theorem}

This will not be the only non-measurable set: in fact if one exists, then there are $2^{2^{\aleph_0}}$ of them, which is of course the maximum amount. Indeed, since the Cantor set (a null set of uncountable cardinality whose existence follows from $ZFC$) has cardinality $2^{\aleph_0}$, and a subset of a null set is null, there will be at least $2^{2^{\aleph_0}}$ different null sets. But adding or subtracting a null set from a non-measurable set doesn't alter its non-measurability, so we'll have $2^{2^{\aleph_0}}$ non-measurable sets. \\ 
\vspace{0,2cm}

Both null and meagre sets are different (and incompatible) accounts of what constitutes a small set of reals. They share a strong structural relationship, known as the measure-category duality, that ensures most arguments are translatable from one to the other. Indeed, by the Erd\H{o}s-Sierpi\'{n}ski Duality Theorem, assuming $CH$ all arguments are translatable, thanks to an involution in $\mathbb{R}$ that swaps null and meagre sets (in fact, assuming the weaker Martin's Axiom suffices) (19 in \cite{Oxt}). So of course we will have this duality in $L$.

But we don't even need this assumption: we always have the category analogue of Fubini's Theorem, which is the Kuratowski-Ulam Theorem, and from that we can reconstruct the argument for category.

\begin{definition}
A \emph{meagre subset of $\mathbb{R}^n$} is one which can be expressed as the countable union of nowhere dense subsets. A \emph{co-meagre set} is one with meagre complement.
\end{definition}

\vspace{-0,3cm}

\begin{theorem} {\bf (Kuratowski-Ulam Theorem for meagre sets)}
 \\ Suppose $W\subseteq \mathbb{R}^2$ has the Baire property, and for $x \in \mathbb{R}$ let  $W_x = \{ y\: | \: (x, y) \in W\}$. \\ Then $W$ is meagre iff $\{ x \in \mathbb{R} \: | \: W_x$ is meagre in $ \mathbb{R}\}$ is co-meagre in $\mathbb{R}$.
\end{theorem}

\begin{corollary} {\bf (1.2 in \cite{Cai})}
Any well-ordering of a non-meagre set of reals does not have the Baire property.
\end{corollary}

Of course, the complexity of the set $W$ in the proof of Theorem 2.3 is that of the formula defining the well-order, so depending upon how easily definable it is we'll get simpler or more complex non-measurable and non-Baire sets of reals. In $L$ these sets are found at the $\Delta^1_2$ level, which is astoundingly low, lower than the complexity of many sets found in everyday mathematical practice. So in accepting $V=L$ some very foundational and far reaching tools for analysis and algebra are at risk of breaking down easily.

As an exemplification, the following sets are $\Pi^1_2$-complete, and thus can't have lower complexity than $\Delta^1_2$ (37 in \cite{Kec}):  \vspace{0,2cm} \\   \indent $\{f \in \mathcal{C}([0, 1]) \: | \: f' \mbox{ is exhaustive}\}$ \\
 \indent $\{f \in \mathcal{C}([0, 1]) \: | \: f \mbox{ satisfies the Mean Value Theorem}\}$ \\
 \indent $\{K \subseteq \mathbb{R}^3 \: | \: K \mbox{ is compact and path connected}\}$
 \\

\section{Sierpi\'{n}ski and Luzin sets}
It might not seem shocking that certain pathological properties can be derived as in the previous section, since a formula well-ordering the universe can be considered very counter-intuitive, conflicting with the assumed vastness and richness of the set theoretic universe. But now we present another odd result derived from a much more modest claim, the Continuum Hypothesis.

$CH$ is equivalent to every infinite subset of real numbers being equinumerous either to the natural numbers or the whole of the real numbers. This assertion seems way more plausible, and in fact the intuitive and intensional understanding of the real line might favour it. After all, why should another cardinality exist in between these? The natural numbers already serve as our archetype for a discrete infinity, and the reals step in as that for a continuous, spatial one. If we take this intuitive understanding of infinities at face value, and especially if we relate them to their use in other mathematical fields, there would seem to be no practical need for another cardinality. Of course, $CH$ remains nonetheless undecided, thanks to $ZFC$'s indeterminacy of the power set function. And in fact, as we will see, its $CH$ being true yields the following set.

\begin{definition}
A \emph{Sierpi\'{n}ski subset of $\mathbb{R}$} is an uncountable set whose intersection with every null set is countable.
\end{definition}

The oddness of this set originates now from a discrepancy between two notions of size different to the previously mentioned: that between cardinality and measure. The Sierpi\'{n}ski set is bijectable with the whole of the reals, and yet manages to coincide with any one of all the null sets in only countably many points.

Of course, this is not trivial because uncountable null sets do exist, such as the Cantor set. In fact, the Cantor set is itself an example of disagreement between size notions (since it is big in cardinality yet null and meagre), and its existence does follow from $ZFC$, so maybe sets showcasing these discrepancies are more fundamental than they might seem. But the Sierpi\'{n}ski set does so in a different way, involving the whole of $\mathbb{R}$ and its subsets, and maybe this difference proves relevant.

\begin{theorem} {\bf (Luzin, 4.3 in $\cite{Sie}$)}
 Assuming $CH$, a Sierpi\'{n}ski set exists.
 \begin{proof}
  Every null set is contained in a null $G_\delta$ set which is a countable intersection of open sets. And since every $G_\delta$ can be coded by a subset of $\omega$ (consider the endpoints of intervals), by $CH$ we can enumerate in order-type $\aleph_1$ all $G_\delta$ sets. For every ordinal $\alpha < \aleph_1$ choose a real $x_\alpha$ not in any one of the $G_\delta$ sets indexed by $\delta < \alpha$, which is possible since their union is null (by being a countable union of nulls), and thus not the whole of $\mathbb{R}$. The uncountable set $X$ of all these reals has only countably many elements in each $G_\delta$, and thus in each null set.
 \end{proof}
\end{theorem}

As in the previous section, a similar argument regarding the Baire property instead of Lebesgue measurability shows that $CH$ implies the existence of a Luzin set, the analogue in category to the Sierpi\'{n}ski set.

\begin{definition}
A \emph{Luzin subset of $\mathbb{R}$} is an uncountable set whose intersection with every meagre set is countable.
\end{definition}
\vspace{-0,8cm}
\begin{corollary}
 Assuming $CH$, a Luzin set exists.
 \end{corollary}

We've seen some examples of how adding axioms to $ZFC$ restricting the richness and vastness of the universe can spawn sets with unexpected properties. $ZFC$ itself proves the existence of a Cantor set, and by assuming further $CH$ (which partially simplifies cardinal arithmetic) we get a Sierpi\'{n}ski set and a Luzin set. $ZFC$ itself also proves the existence of a non-measurable, non-Baire set, but the more structural regularity there is in our universe (and the easier it is to define this regularity), the lower its complexity.

\chapter{Other inner models}

G\"{o}del's strategy for constructing the smallest model of $ZF$ was promptly extended to allow for models containing more sets and satisfying stronger theories. This provided a rich spectrum of inner models with ever increasing complexity. And analogously of different Axioms of Relative Constructibility with ever increasing consistency strength.

Of course, all of these models could only be built inside an already existing set theoretic universe, and so building them is usually tantamount to finding a smaller universe with more regularity and less complexity than the one we started with. Developments in the complementary direction had to wait for Cohen's revolutionary technique of Forcing in the 60s, which through a far less direct construction allowed for increasing the universe's complexity. These two approaches thus provide thorough tools for the study of possible set theoretic universes.

But there's another remarkable use for Inner Model Theory. The theories of inner models present regularities which allow for finer analysis, elucidating many otherwise intractable issues, in a way similar (but more complex) to $L$'s Fine Structure. So being able to conceptualize a model (of a certain strong theory) as actually an inner model of another, bigger universe will help answer some questions. That's why the search for canonical inner models of ever stronger large cardinal axioms is a central program to modern set theory, which includes Woodin's program searching for Ultimate-$L$, a canonical inner model where a supercompact cardinal exists \cite{Sar} \cite{Woo}. G\"{o}del had already remarked upon the role such axioms of infinity could play in future developments ($§$3 in \cite{Kan}).

In this Chapter we only consider inner models of relatively low consistency strength inside the large cardinal hierarchy of axioms, which can be built in a straight forward way. We'll now see there are two different manners of constructing an inner model from a set or class by paralleling the construction of $L$.

\vspace{0,1cm}

\begin{definition} 

 \phantom{a} \\ $\mathcal{D}(B)$ is the set of all subsets of B definable with parameters in B. \\That is, of all $C\subseteq B$ such that, for a certain formula $\phi$, $C = \{ x \in B \: | \: B \models \phi(x) \}$.
 
 \vspace{0,1cm}
 
 \noindent $\mathcal{D}^A(B)$ is the set of all subsets of B definable in B with parameters, and additionally using a predicate added to the language which is interpreted as $A$.\\ That is, $\phi$ can now include the predicate $A(x)$, interpreted as $x \in A$.
\end{definition}

Notice this definition is not equivalent to that of $Def^A(P)$ in Chapter 1. There we were defining elements, and here subsets. So that $Def^A(P) \subseteq A$, while $\mathcal{D}^A(B) \subseteq \mathcal{P}(B)$. From it we can define the following hierarchies of relative constructibility:

\begin{definition}
 
For $A$ a set (or also a proper class in the left column) \begin{multicols}{2}
$L_0[A]=\varnothing$

$L_{\alpha+1}[A]=\mathcal{D}^A(L_{\alpha}[A])$

$L_{\gamma}[A]=\bigcup_{\alpha<\gamma}L_{\alpha}[A]$ for limit $\gamma$

$L[A]=\bigcup_{\alpha}L_{\alpha}[A]$

$L_0(A)=TC(\{A\})$

$L_{\alpha+1}(A)=\mathcal{D}(L_{\alpha}(A))$

$L_{\gamma}(A)=\bigcup_{\alpha<\gamma}L_{\alpha}(A)$ for limit $\gamma$

$L(A)=\bigcup_{\alpha}L_{\alpha}(A)$

\end{multicols}

\end{definition}
 
 That is, $L(A)$ just adds all of the elements and information in $A$ at the beginning, while $L[A]$ only uses $A$ in every level for additional definability power, so that maybe $A \nsubseteq L[A]$\footnote{This does sometimes happen, for instance $\mathbb{R} \nsubseteq L[\mathbb{R}]$ when a Cohen real over $L$ exists \cite{Coh}.}. But both will always be an inner model of $ZF$ by a proof analogous to that of $L$ (II.7.2B in \cite{Dev}, II.6.30 in \cite{Kun}). $L[A]$ is usually interpreted model theoretically as the structure $\langle L[A], \in, A \rangle$, with the additional predicate used in building it also present.
 
 Notice in the case of $L(A)$, the transitive closure is done on $\{A\}$ to have $A \in L(A)$. This ensures any definition carried out in $L[A]$ can also be performed in $L(A)$ (by using the parameter $A$ instead of the predicate $A(x)$), and so clearly $L[A] \subseteq L(A)$. Notice also that $TC(\{A\})$ is not a set for $A$ a proper class, and thus the construction isn't possible then.
 
 In Inner Model Theory the mainly used construction is $L[A]$, since its minimal definability properties ensure it will be the smallest inner model satisfying a certain theory, unlike $L(A)$. It also allows for a stronger analogue to $L$ in results like the following:
 \vspace{-0,3cm}
 \begin{theorem}{\bf  (Generalized Condensation) (II.7.4A in  \cite{Dev})}
\label{bottLine}
If $\alpha$ is a limit ordinal, $\pi$ is the Mostowski collapse and $X \prec_1 L_\alpha[A]$, then there is a unique $\beta\leqslant\alpha$ such that
\begin{itemize}
    \item[] $\pi: \langle X, \in\rangle \cong \langle L_\beta[\pi(A)], \in\rangle$
\end{itemize}
So in particular if $A$ is transitive, we recover regular Condensation.
\end{theorem}

 \begin{theorem}{\bf  (Generalized Partial GCH) (II.7.4C-E in  \cite{Dev})}
\label{bottLine}
\\ If $V=L[A]$ and $A \subseteq \kappa^+$, then $\forall \lambda \geqslant \kappa (2^\lambda = \lambda ^+)$. \end{theorem}

This last result is especially useful because it can be seen that every $L[A]$ with $A$ a set is the same model as another $L[A']$ where $A'$ is a set of ordinals (see the proof of 3.10.5 for an example).
 
 We proceed to the motivation and description of concrete inner models.

\section{L[\#$_{1}$]}

Motivated by model-theoretic results by Ehrenfeucht and Mostowski \cite{Ehr}, the study of $L$-indiscernibles culminated in the isolation of $0^\#$, a crucially canonical object whose existence presents sweeping consequences for the set-theoretical universe. $0^\#$ is, broadly speaking, a set of formulas (coded as a real through G{\"o}del numbering) coding the theory of $L$. It of course doesn't belong to $L$, and in fact its existence (which is independent of $ZFC$) entails $V$ to be vastly larger than $L$ in many relevant aspects, and conversely its non-existence entails $V$ and $L$ to be way more similar.

An indepth exposition of $L$-indiscernibles and $0^\#$ can be found in V of \cite{Dev} or 9 of \cite{Kan}. We won't summarize here that exposition due to lack of space, but the three following results might exemplify the role of $0^\#$:

\vspace{-0,3cm}

\begin{lemma}
{\bf (9.17 in  \cite{Kan})}
If $0^\#$ exists, \\ then $|\mathcal{P}(x)^L| = |x|$ for every infinite $x\in L$, so in particular $\mathcal{P}(\omega)^L$ is countable.
\end{lemma}

\vspace{-0,3cm}

 \begin{theorem} {\bf (Kunen, see V.4 in  \cite{Dev})}
  \\ $0^\#$ exists iff there is a nontrivial elementary embedding L $\prec$ L
\end{theorem}

\vspace{-0,3cm}

 \begin{theorem} {\bf (Jensen's Covering Theorem, V.5.1 in  \cite{Dev})}
  $0^\#$ does not exist iff \\ for any uncountable subset $X \subseteq ON$, there is a $Y\in L$ with $Y\supseteq X$ and $|Y| = |X|$

\end{theorem}

  \vspace{0,5cm}
 
 This treatment can readily be generalized to a study of $L[x]$-indiscernibles and the set of formulas $x^\#$, for any set $x$ (although results like the previous ones aren't completely translated). We can adjoin some of these sharps to $L$ to obtain slightly larger models, the first of which is $L[0^\#]$. These models will generally have a really similar structure to $L$, as do all of the $L[A]$ for $A\subseteq\omega$.

The next model which presents a considerable step up is $L[\#_1]$, the smallest model closed under the sharps of its reals. It serves as a canonical framework of richer reals, and as we'll see the higher complexity helps alleviate some pathologies.
 
 So we'll be particularly interested in the sharps of reals. By defining $x^\#$ as an E-M set for $L[x]$ (see $§$14 in \cite{Kan}), "$x^\#$ exists and $x^\#=y$" becomes a $\Pi^1_2$ formula without parameters, that is, of the form $\exists z \forall t \phi(x, y, z, t)$, where the quantifiers range over real numbers (14.16 in \cite{Kan}). Thanks to Shoenfield's Lemma, these sharps are absolute for inner models, in the following sense.
 
   \begin{lemma} {\bf (Shoenfield's Absoluteness Lemma, \cite{Sho})}
  \\ Any two inner models agree on the truth of $\Pi^1_2$ sentences.
\end{lemma}

 \vspace{-0,3cm}
 
 \begin{corollary}
For $x, y \in M$, \:\:\: $M \models$ "$x^\#$ exists and $x^\#=y$"  \:\:\: iff \:\:\: $x^\#=y$ \end{corollary}

  \vspace{0,5cm}
  
 We're interested in studying the smallest inner model (and thus a model of the form $L[A]$) closed under real sharps. That is, the smallest inner model satisfying $L[A] \models \forall x \subseteq \omega(x^\#$ exists$)$. As we'll see, this model might or might not exist depending on how many real sharps exist in our set theoretic universe.
 
 But let us remark that the existence of a measurable cardinal (a relatively modest large cardinal) does imply the existence of all real sharps \cite{Kan}, and so under that assumption, or in any model with a measurable cardinal, the construction presented below will yield the desired model.
 
 \vspace{0,4cm}
 
 Consider first the sharp function on the reals, $F_1$.\footnote{The 1 stands for the reals being \textit{objects of class one} \cite{Dub}, as usually we talk about the naturals being class zero, the reals class one, etc.} That is,
\begin{align*}
F_1 \colon \mathcal{P}(\omega) &\to \mathcal{P}(\omega) \\
x &\mapsto x^\#
\end{align*}

 This is a partial function on the reals, and is only a total function when all of the sharps of reals exist. Now, the obvious construction with our previous definition of $L[A]$ won't work, because $L[F_1] = L$. Indeed, $\mathcal{D}^{F_1}(L_\alpha[F_1]) = \mathcal{D}(L_\alpha[F_1])$, since  $$L_\alpha[F_1] \models a\in F_1 \:\:\: \textnormal{iff} \:\:\: L_\alpha[F_1] \models a=\langle x, y \rangle \wedge x^\#=y$$
 which was already a definable predicate as explained above.
 
 On the other hand, $L(F_1)$ contains all of the sharps of reals existing in our universe, and thus might not be the \textit{smallest} model closed under real sharps as required\footnote{This does indeed sometimes happen. If $c$ is a Cohen real over $L(F_1)$ (\cite{Coh}), in $L(F_1)[c]$ still every real has a sharp, but also in $L(F_1)[c]$ the smallest model closed under real sharps doesn't contain $c$.}.
 
 So consider instead, as in \cite{Bag} (but with different notation), \begin{itemize}
     \item [] \: \: \: \:\: \:\: \:\: \:\: \:\: \:\: \: $\#_1 = \{ \langle x, n \rangle \in \mathcal{P}(\omega) \times \omega \, | \, x^\#$ exists $\wedge \: n \in x^\#\}$
 \end{itemize}
 
 Then $L[\#_1]$ is the appropriate model, because we've made the information encoded in the sharps available to use for successive definitions.
 
 $L[F_1]$ could also work by changing its definition slightly. As in \cite{Dub}, define  $\hat{L}[F_1]$ by changing the second clause in Definition 3.2 to
 $$\hat{L}_{\alpha+1}[F_1]= \mathcal{D}(\hat{L}_{\alpha}[F_1]\cup \{x^\# = F_1(x) \: | \: x \in \hat{L}_{\alpha}[\#_1]\cap \textnormal{Dom}(F_1)\})$$
 
 That is, at each stage we explicitly add the sharps of all of the already constructed reals (something $L[\#_1]$ already did automatically).
 
 We present some useful basic facts about these models. The above results 3.3 and 3.4 are especially useful when $A$ doesn't have good definability properties. But in the present case the definition of $\#_1$ already provides a direct proof of Condensation.

 \begin{lemma} 
  \begin{enumerate}
     \item $\hat{L}[F_1] = L[\#_1]$
     \item If M is an inner model satisfying M $\models \forall x \subseteq \omega(x^\#$ exists), then $L[\#_1] \subseteq M$
     \item $L[\#_1]^{L[\#_1]} = L[\#_1]$
     \item Condensation is valid for $L[\#_1]$
     \item $L[\#_1] \models GCH$
     \item There's a global well-ordering of $L[\#_1]$ $\Delta^1_3$-definable in $L[\#_1]$
     \item For $a \subseteq \omega$, if $a^\# \in L[\#_1]$ then $a \in L[\#_1]$
     \item $L(F_1) = L(\#_1)$
     \item If M is an inner model satisfying $\forall x \subseteq \omega(x^\#$ exists $\rightarrow x^\# \in$ M), then $L(\#_1) \subseteq M$
 \end{enumerate}
 
\begin{proof}
 \begin{enumerate}
     \item If $a \in L_\alpha[\#_1]\cap\mathcal{P}(\omega)$, and $a^\#$ exists, then $a^\# = \{n \in L_\alpha[\#_1] \: | \: L_\alpha[\#_1] \models \langle a, n \rangle \in \#_1 \} \in L_{\alpha+1}[\#_1]$. So by induction on $\alpha$, $\hat{L}[F_1] \subseteq L[\#_1]$. For the other inclusion, suppose $L_\alpha[\#_1] \in \hat{L}[F_1]$ and $a \in L_{\alpha+1}[\#_1]$. Then \begin{itemize}
     \item[] $a = \{x \in L_\alpha[\#_1] \: | \: L_\alpha[\#_1] \models \phi (\bar{y}, \#_1, x) \} =$
     \item[] $= \{x \in \hat{L}[F_1] \: | \: \hat{L}[F_1] \models (L_\alpha[\#_1] \models \phi (\bar{y}, \#_1, x))\}$
     \end{itemize}
     
     Since $\#_1\cap L_\alpha[\#_1] \in \hat{L}[F_1]$, and given the absoluteness of the definition of $L_\alpha[\#_1]$ and the satisfiability predicate, by Comprehension this set belongs to $\hat{L}[F_1]$.
     \item We see by induction on $\alpha$ that $L_\alpha[\#_1]^M = L_\alpha[\#_1]$. If it's true for $\alpha$ and $a \in L_{\alpha+1}[\#_1]$, then $a = \{x \in L_\alpha[\#_1]^M \: | \: L_\alpha[\#_1]^M \models \phi (\bar{y}, \#_1, x) \}$, so
     \begin{itemize}
     \item[] $M \models a = \{x \in L_\alpha[\#_1] \: | \: L_\alpha[\#_1] \models \phi (\bar{y}, \#_1, x) \}$
     \end{itemize}
     because $\#_1\cap L_\alpha[\#_1]^M \in M$, so $a \in L_\alpha[\#_1]^M$.

     \item One inclusion is immediate, for the other we just use the absoluteness with respect to inner models of the satisfiability relation and of the definition of $\#_1$. This is assured by the definition of $\#_1$ being $\Pi^1_2$, and Shoenfield's Absoluteness Lemma. \\
     In fact, when $L[\#_1]$ is indeed closed under real sharps, 3 is also directly implied by 2, since $L[\#_1]^{L[\#_1]}$ is also an inner model of $V$.
     
     \item By recreating the original proof of Condensation (as in our source II.5.2 of \cite{Dev}), the only crucial change is noticing there is a $\Sigma_0$ formula $\phi$ such that \\ $v = L_\gamma[\#_1] \leftrightarrow \exists z \phi(z, v, \gamma)$. Indeed, the definition of the hierarchy is obtained by just implementing that of $\#_1$, and all quantification can still be bounded by a certain $z$ since we're only dealing with members of $\mathcal{P}(\omega)$ and their collections.

     \item Just as in the source for Lemma 1.6, the hierarchy satisfying Condensation implies it proves $GCH$. Equivalently, it is because (as for $v = L_\gamma[\#_1]$ in 4) we can write the equivalent axiom of relative constructibility $V=L[\#_1]$ and that proof of $GCH$ goes through (even if the complexity has been raised as we'll see next).
     
     Alternatively, and as an illustrative example, by Theorem 3.4 we only needed to see $\#_1$ can be coded as a subset of $\omega_1$. Indeed, given an inclusion $i$ of $\omega$ into $\omega_1$ and a pairing function $J_1$ for $\omega_1$ (both constructible), we can define $A = \{ J_1(x, i(n)) \: | \: x^\#$ exists $\wedge \: n \in x^\#\}$.
     
     \item Exactly as for $L$, with $\#_1$ now an additional symbol of the language. Again, by keeping track of the complexity through the construction, we can see the definition of $\#_1$ being $\Pi^1_2$ now makes the well-ordering $\Delta^1_3$ (check the sources for Theorem 1.4).
     \\ The well-ordering in $L_\alpha[\#_1]$ can again be made definable over $L_\alpha[\#_1]$ by the construction of Boolos.
     
     \item If G is a G{\"o}del numbering, $a = \{n\in\omega \: | \: G(\floor{n}\in a) \in a^\# \}$, where $\floor{n}$ is the term without free variables representing the numeral $n$, and thus $\floor{n}\in a$ is a sentence of $L[a]$. Since G and $\floor{\,}$ are definable functions, $a \in L[\#_1]$.
     
     \item We just need to show each object belongs to the other inner model. We build over the $\omega_2$ levels to ensure all necessary sets are present:
     
     $\#_1$ = \{$a \in L_{\omega_2}(F_1) \: | \: L_{\omega_2}(F_1) \models a = \langle x, n \rangle \wedge x \in $Dom$(F_1) \wedge n \in F_1(x)$\} $\in L(F_1)$ 
    
     $F_1 = \{ a \in L_{\omega_2}(\#_1) \: | \: L_{\omega_2}(\#_1) \models \\ \phantom{aaaaaa}a = \langle x, x' \rangle \wedge \forall n\in\omega(n\in x' \leftrightarrow \langle x, n \rangle \in \#_1) \wedge x' \neq \emptyset \} \in L(\#_1)$
     
     \item  Since $M$ contains all of the real sharps, and is a transitive model of $ZF$, by Replacement $TC(\{\#_1\}) \in M$. By the absoluteness of the $\mathcal{D}$ function on sets and the definition of $\#_1$, and by induction on $\alpha$, $L_\alpha (\#_1) ^ M$ = $L_\alpha (\#_1)$.
     \end{enumerate} \end{proof}
 \end{lemma}
 
Another interpretation of $L[\#_1]^{L[\#_1]} = L[\#_1]$ is that, given that there exist enough real sharps for $L[\#_1]$ to be closed under them, the model will remain the same no matter how many more real sharps actually do exist. So we might as well study it with the assumption that all real sharps exist.

\vspace{0,2cm}

An interesting question is how many (and which) real sharps actually have to exist for $L[\#_1]$ to be closed under them. This is equivalent to asking how many reals are in $L[\#_1]$ when it's thus closed, and that's equivalent to determining $\omega_1^{L[\#_1]}$, since by the $GCH$ $|\mathcal{P}(\omega)|^{L[\#_1]} = (2^{\aleph_0})^{L[\#_1]} = \aleph_1^{L[\#_1]} $. This will of course be bound by $\aleph_1^{L} \leqslant \aleph_1^{L[\#_1]} \leqslant \aleph_1$. \\
\indent Now, if there are enough real sharps for $L[\#_1]$ to be closed under them, then since $L[\#_1]^{L[\#_1]} = L[\#_1]$ it is a model for $\aleph_1^{L[\#_1]} = \aleph_1$, so this equality can't be refuted in $ZFC + \forall x \subseteq \omega(x^\#$ exists$)$.  \\
\indent But we can say more: it is independent of $ZFC + \forall x \subseteq \omega(x^\#$ exists$)$. Indeed, we can construct a forcing extension satisfying $\aleph_1^{L[\#_1]} < \aleph_1$. Now, in a universe with an unbounded class of measurable cardinals (and so assuming the consistency of this statement) necessarily every sharp exists \cite{Kan}. Without delving into the complex technique of Forcing, suffice it to say that then we can apply Forcing to collapse $\aleph_1^{L[\#_1]}$ to a countable ordinal, while maintaining the existence of the unbounded class of measurables (since the Forcing will only alter an initial segment of the ordinals), and thus the existence of every sharp.

\vspace{0,2 cm}

Adding to $L$ an amount of sharps smaller than $\aleph_1^{L[\#_1]}$ (let alone $\aleph_1^{L}$) will clearly not suffice to close it under real sharps. But notice why this happens: adding for instance all of the successive 0$^{\#n}$ doesn't suffice, since then we'd be able to code all of these sharps in $N = \{ \langle n, m \rangle \in \omega \times \omega \, | \, m \in 0^{\#n}\}$, and the model could be expressed as $L[N]$. Since $L[\#_1]$ knows $N$ to be countable, $L[N]$ would too if it were closed under real sharps (since it would contain $L[\#_1]$). But then it would know $N$ is (encodable in) a real, and $N^\# \notin L[N]$.

\section{L[\#]}

We've considered the previous model because our main focus is on the reals, but as seen in Chapter 1, we're really studying the power set operation, and not only the properties of $\mathcal{P}(\omega)$ are relevant to inner models, but also those of $\mathcal{P}^\gamma(\omega)$ (or equivalently $\mathcal{P}(\kappa)$ for other cardinals). In order for these to get the same treatment, it might seem arbitrary for the sharps to be restricted to the reals (sets of natural numbers), and thus we also want to consider the whole sharp function on sets of arbitrary ordinals:
\begin{align*}
F \colon \mathcal{P}(ON) &\to \mathcal{P}(ON) \\
x &\mapsto x^\#
\end{align*}

Notice the image of the function is no longer contained in $\mathcal{P}(\omega)$ (or any $\mathcal{P}(\kappa)$) if enough sharps exist, since ever more parameters will be needed in the resulting theory of $L[A]$, and so the true sentences will have to be codified by a bigger cardinal.

As before, we'll have to write this as a relation for it to add definability power.

\begin{itemize}
     \item [] \: \: \: \:\: \:\: \:\: \:\: \:\:  $\# = \{ \langle x, \alpha \rangle \in \mathcal{P}(ON) \times ON \, | \, x^\#$ exists $\wedge \: \alpha \in x^\#$\}
 \end{itemize}

This relation might now be a proper class, just like $F$.

As mentioned after 3.4, every sharp of a set is actually the sharp of a set of ordinals, and so every sharp of a set will be in the image of $F$. By the same reason, $L[\#]$ will be closed under all sharps.

As before, maybe not all existing sharps are needed to close $L[\#]$ under them, and so maybe $\#\cap L[\#] \subsetneq \#$. For this reason, the following function and relation defined by recursion, as presented in \cite{Wel}, are even more natural:
\begin{align*}
F' \colon ON &\to \mathcal{P}(ON) \\
\alpha &\mapsto (F'\restriction \alpha)^\#
\end{align*}

\begin{itemize}
     \item [] $\#' = \{ \langle \alpha, \beta \rangle \in ON \times ON \, | \, (\#' \restriction \alpha)^\#$ exists $\wedge \: \beta \in (\#' \restriction \alpha)^\#$\}
 \end{itemize}
 
 That is, we're not just taking the sharps of ordinals: we transfinitely iterate the sharp function, starting from 0. At each stage, we take the sharp of everything coming before, so that for instance $F'(\omega) = \{ 0^\#, 0^{\#\#}, \ldots\} ^\#$. This is the sequence of the set theoretical objects called \textit{mice}, which play a crucial role in modern Inner Model and Core Model Theory.
 
 Of course, if not all sharps exist this function and relation will end at some ordinal, and thus will be bounded and not total. Furthermore, Dom($F'$) will be an ordinal (it will be the first $\alpha$ for which $F'(\alpha)$ doesn't exist), since the existence of a sharp implies the existence of all lower sharps.
 
 As mentioned, this recursive definition is more natural because it does ensure $\#'\cap L[\#'] = \#'$. Indeed, any sharp arrived at by the function (through transfinite iteration) will also be arrived at by $L[\#']$, since it also provides transfinitely many levels for definition.
 
 Before proceeding, let us notice that these functions and relations are now $\Pi_2$ instead of $\Pi^1_2$, basically because the $x$ of which we take the sharp can  no longer be bounded as a real.\footnote{More concretely, by going back to the proof of $\#_1$ being $\Pi^1_2$ thanks to definition through E-M sets (14.11, 14.16 in \cite{Kan}), we only need to notice that: \begin{enumerate}
     \item the G\"{o}del numbers can now be coded as members of $\mathcal{P}^\gamma(\omega)$ for a certain $\gamma$ depending on the element $x$ of which we take the sharp, and
     \item the well-orderings ($E_y$) representing the order type of the set of indiscernibles will now possibly have arbitrarily high cardinality (this will be so when $F'$ is total), so we can't assure they're coded by members of $\mathcal{P}^n(\omega)$ for any $n$, and can only be bounded by an arbitrary set, thus dropping the superscript in the complexity class.
     
      \end{enumerate}}
 
 We present as before some basic facts about this model. Most are obtained by reasoning analogous to that of the previous section, so we provide only the proofs with non-trivial changes.
 
 \begin{lemma} 
  \begin{enumerate}
     \item $L[\#] = L[\#']$
     \item If M is an inner model satisfying M $\models \forall x (x^\#$ exists), then $L[\#] \subseteq M$
     \item $L[\#]^{L[\#]} = L[\#]$
     \item Condensation is valid for $L[\#]$
     \item $L[\#] \models GCH$
     \item There's a global well-ordering of $L[\#]$ $\Delta_3$-definable in $L[\#]$
 \end{enumerate}
 
 \begin{proof}
 \begin{enumerate}
     \item [3.] Now we can't use Shoenfield's Absoluteness Lemma since the definition of $\#$ won't necessarily be $\Pi^1_2$ if enough sharps exist. But we can use the other approach mentioned in 3.10.3:
     \\ When L[\#] is indeed closed under all sharps, 3 is implied by 2, since $L[\#]^{L[\#]}$ is an inner model of $V$. And even when it isn't, since $F'$ will have a univocal cut-off point $\alpha$, we will have 2 for the formula "Dom($F'$) = $\alpha$" (using $\alpha$ as a parameter), and this will also imply 3.
     \item [4.] As before, but noticing now why all quantification can be bounded by a certain $z$: because in the level $L_\gamma[\#]$ the highest sharp possibly appearing is $F'(\gamma)$, and so we'll be dealing with members of at most $\mathcal{P}^{\gamma'}(\omega)$ for a certain $\gamma'$, which can be bounded by a set.
     \item [5.] Of course, Theorem 3.4 can no longer be used because we're not assured the relation will be coded as a subset of $\omega_1$ (or even of any cardinal). But the usual proof through the formula $V = L[\#]$ still applies.
 \end{enumerate}
 \end{proof}
 \end{lemma}

\section{Gaps in inner models}

The generalization of gap results for these two cases is immediate: all developments in Chapter 1 can be applied straight away to $L[\#_1]$ and $L[\#]$. This is because they both satisfy all of the good properties needed as tools for the study of gaps: mainly Condensation (and thus also $GCH$) and the global well-order (and thus also $AC$).

As seen above, these are consequences of $\#_1$ and $\#$ being univocally definable without parameters, so that we can write the axiom $V=L[A]$ (or the formulas $v=L_\gamma[A]$) just as in the case for $L$. This also ensures the absoluteness between inner models of the notion of "being $L_\alpha[A]$", and the validity of the Skolem hull arguments used (as in 1.19). 

So it seems like it is this definability what binds these models so close to $L$ regarding gaps. Even though the existence of sharps does provide richer real numbers (or elements of $\mathcal{P}^\gamma(\omega)$), if they do exist they are nicely definable. And this gives the inner models the regularity properties required.

As a more concrete exemplification, we can intuitively see why adding the sharp relations won't alter the gap structure. They will only contribute to the hierarchy by adding a real's sharp to the level immediately after that real appears (as especially elucidated by the definition of the first model as $\hat{L}[F_1]$). But this will not stop any gap from appearing. Indeed, if a sharp is added by the function at a level, then the previous level added a real. And since successor levels can't start gaps, this means that level wouldn't be a gap level anyway.

\vspace{0,4cm}

What ensures the results provable is both the expressibility power of the reals in the model on the one hand, and the regularity properties of the model on the other. For inner models bigger than $L$, the expressibility required will never fail. Indeed, the definability power over the levels of the model will be equal or greater (augmented by a predicate), so we can always codify information in the reals as before. So only the regularity properties can fail.

It is therefore natural to wonder when they fail. That is, when aren't the levels of the hierarchy so neatly definable? Or how big must the model be for, even more strongly, Condensation or $GCH$ to fail?

Apparently not very big: as exposed in \cite{Wel}, an inner model with an $\omega_1$-Erd\H{o}s Cardinal won't satisfy Condensation. The existence of these cardinals is just slightly stronger than $\forall x  (x^\#$ exists$)$, and weaker than the existence of a Ramsey cardinal, so many inner models do satisfy it, just above the two we've studied.

For instance, since the existence of a measurable cardinal is even stronger, the minimal model $L[U]$ for a measurable cardinal doesn't satisfy Condensation, so the gap results won't be generally applicable (although it does satisfy $GCH$ below its measurable cardinal, and thus weaker results there might be possible \cite{Mit}).

\vspace{0,2cm}

Another result of \cite{Wel} is even more interesting for our purposes: a model satisfying $GCH$ is equivalent to another property essential for our study of gaps.

\begin{definition}
 {\bf (Acceptability)}
 For $A \subseteq ON$, we say the model $L[A]$ is acceptable if, whenever $B \in \mathcal{D}^A(L_\alpha[A])\cap\mathcal{P}(\rho)$ and $B\notin L_\alpha[A]$, \\ we have $\exists F \in L_{\alpha+1} (F\colon \rho \rightarrow L_\alpha[A]$ is an onto function$ )$
\end{definition}
\vspace{-0,6cm}
\begin{theorem}
For $A \subseteq ON$, $L[A]$ is acceptable iff $L[A] \models GCH$
\end{theorem}

Notice acceptability is very reminiscent of our notion of $\beta$-analytical copies, and the isomorphisms witnessing them. Acceptability states that whenever a new subset of $\rho$ appears in a level (analogously, when $\alpha$ is a gap of a certain order) then this level knows the previous level to be injectible into $\rho$ (that is, of cardinality at most that of $\rho$). The definition doesn't require this injection to be an isomorphism, as was our case. But of course, if acceptability fails, then the stronger theorem by Boolos demanding it be an isomorphism (our Lemma 1.19) will also fail. And so we are left with no apparent way of proving the start of a gap must be a limit. This completely disrupts the gap regularity we had found in Chapter 1.

So a failure of $GCH$ makes the study of gaps as we know it impossible. As seen above, this can only happen if Condensation fails. And as mentioned before, even if only Condensation fails, results as basic as the existence of arbitrarily big gaps below a cardinal (our Theorem 1.9) would need a completely different proof, if one does exist.

This is a general pattern which also affects deeper results about inner models by Jensen. In the words of Welch, without some form of Condensation, fine structural analysis is hopeless. \newpage

\section{Pathologies in inner models}

We see now how the two results of Chapter 2 generalize to these and other inner models.

\vspace{0,2cm}

Since $GCH$ still holds, a Sierpi\'{n}ski set and a Luzin set will of course still exist in both $L[\#_1]$ and $L[\#_1]$. But there are many inner models in which this can fail: according to Theorem 3.4, choosing $A$ containing uncountable ordinals might suffice. But it still might be that $L[A]$ possesses some good properties (like definability in our case) that ensure $GCH$. And indeed, even big models like $L[U]$ satisfy $CH$.

In fact, $\neg CH$ is a property that seems more natural for a Forcing extension than an inner model (and indeed, Forcing extensions were invented for that in the first place), since through Forcing we can collapse cardinals and thus alter cardinal arithmetic. This is augmented by the fact that $ZFC$ alone doesn't even prove there exists an inner model satisfying $\neg CH$, since this won't happen if $V=L$, which is consistent with $ZFC$. Furthermore, some of the central regularity properties characteristic of many inner models imply $CH$, like for instance diamond ($\diamond$). Even more, $AD$ being true in the universe implies all inner models of a certain natural form satisfy $CH$ \cite{Bec}.

Nonetheless, inner models falsifying $CH$ are in general possible. For instance, if $CH$ is false in the universe, then not only is $V$ trivially an inner model of $\neg CH$, but also $L(\mathbb{R})$.\footnote{Or more rigorously, of a formulation of $\neg CH$ which employs subsets of reals instead of cardinal arithmetic, since $AC$ might fail in $L(\mathbb{R})$ and thus we can't talk about $(2^{\aleph_0})^{L(\mathbb{R})}$.}

When $CH$ thus fails, the proof provided won't work, and thus it might be a priori that a Sierpi\'{n}ski or Luzin set doesn't exist. Still, to ensure it doesn't, we need strong large cardinal axioms. For instance, Martin's Axiom for $\aleph_1$ does imply a Sierpi\'{n}ski set can't exist (and thus also $\neg CH$) (V.6.29 in \cite{Kun}), but this again is a statement that holds in Forcing extensions.

All in all, the existence of such sets is pervasive amongst inner models. Given, as remarked earlier, that its counter-intuitive features are not much greater to those of a Cantor set (which exists in $ZFC$), this might lead us in the direction of considering them not as pathological as they might have seemed. \\

 Regarding non-measurable and non-Baire sets, both $L[\#_1]$ and $L[\#]$ do of course still have one: the well-ordering of their reals. What changes is the complexity of their definition. This change is only possible because there's not a single formula univocally defining a well-ordering of the reals in the theory $ZFC$ (that is, in every model of $ZFC$).
 
 As seen above, in $L[\#_1]$ these sets will be at most $\Delta^1_3$ definable, while in $L[\#]$ we can only assure they are $\Delta_3$. This difference is very relevant. A $\Delta^1_3$ set is still a pretty natural construction in second-order arithmetic, and thus the situation is just marginally less pathological than that of $L$. A $\Delta_3$ set on the other hand possibly can't be defined in any arithmetic, and indeed its quantifications might range over arbitrarily big ordinals. This is thus a construction involving the full vastness of the set theoretic universe.
 
 Of course, this a priori doesn't exclude the possibility that some well-order on the reals, defined in a way different from the usual proof, does have lower complexity, let alone any non-mesurable or non-Baire set whatsoever.
 
 But for the case of the well-order we do have relevant limitative results, which furthermore are related to the structure of the models (the complexity of the well-orderings of the reals is an important issue in Descriptive Set Theory, and much intersected with Inner Model Theory). The first of them is the following:
 
 \vspace{-0,3cm}
 
 \begin{theorem}  {\bf (Mansfield, 25.39 in  \cite{Jec})}
 If for a certain $A$ there is a $\Sigma^1_2(A)$ well-ordering of the reals, then every real is constructible from $A$, that is, belongs to $L[A]$.
 \end{theorem}
 
 This implies that the $\Delta^1_2$ well-ordering of the reals is only present in $L$. In fact, this complexity is only possible because $\mathbb{R}^L$ itself is $\Delta^1_2$ in $L$. As we'll see shortly, this is a repeating pattern: in many inner models the complexity of the well-ordering is exactly that of the set of reals itself. So in a sense it is the complexity of $\mathbb{R}$ itself what keeps the well-ordering (and thus the non-measurable set) from being too simple.
 
 This phenomenon is due to these inner models being \textit{canonical}, in the sense that they are definable without parameters, that is, $M = \{ x \: | \: \phi(x)\}$ for $\phi$ without parameters (for instance, $\phi(x) \equiv x\in L$). This definition grants them a certain locality: they can be built inside any model. This includes themselves, and so they can reconstruct their building from the inside, yielding a well-order. That's why the well-order will have the same complexity as the constructed set itself: because the construction of the set automatically yields a well-order.
 
 \vspace{0,3cm}
 
 As an application of Mansfield's result, the existence of a measurable cardinal implies that there is no $\Delta^1_2$ well-ordering. But the canonical minimal model $L[U]$ for a measurable cardinal does have a well-ordering of its reals of complexity $\Delta^1_3$. In a similar spirit, if a model contains a Woodin cardinal (a relatively strong large cardinal axiom), then its $\mathbb{R}$ can't be $\Delta^1_3$ definable, and so neither its well-order. By adding further Woodin cardinals we keep rising this complexity, until at infinitely many Woodin cardinals there is no projective well-ordering of the reals whatsoever.
 
  \vspace{0,3cm}
  
  By considering non-canonical inner models we can observe different behaviour. We can even consider a model where there is no well-ordering of the reals at all (and so necessarily Choice fails).
  
  For instance, assuming $AD$ (incompatible with $AC$), not only does $L(\mathbb{R})$ not have a well-ordering of the reals (we have lost it when adding $\mathbb{R}$ at the beginning of the construction): it doesn't have any non-measurable sets at all.
  This model was brought into the spotlight because Kechris showed that, under the assumption of a supercompact cardinal, it satisfies $AD$ and also Dependent Choice, a weaker version of Choice. These are in some sense the best analytical properties we can hope for in a model without non-mesurable sets, so as Jensen put it this is an analyst's dream.
  
  But of course, on the other hand, a non-canonical model won't have the minimality properties to be as fundamentally relevant to the development of set theory as the models we've mainly been dealing with.
  
   \vspace{0,3cm}
   
   In conclusion we see that, even if not all inner models do, canonical inner models present a direct correlation between the complexity and vastness of the model (and the consistency strength of its theory) and the definitional complexity of its reals and their well-ordering.

\newpage
\thispagestyle{empty}
\addcontentsline{toc}{chapter}{Conclusions}
\chapter*{Conclusions}

The results of Chapter 1 demonstrate that a certain simpler form of fine structural analysis is also possible and fruitful for the G\"{o}del hierarchy, and well suited to prove some results about the power set operation. Generalized gaps can be applied in many more directions than the ones presented, and they can be a useful tool for inner models.

The reals, thanks to their multipurpose application in definitions (due to their canonical definition), are deeply linked with the structure and construction of the model they inhabit. More generally, the complexity of the power set function in an inner model very accurately gauges the richness of said model.

In $L$ the regularity is utmost and observable in any power set. Regarding $L[\#_1]$ and $L[\#]$, as much as sharps might present a richer paradigm for model theoretic study, they're not enough to disrupt the regularity, since they are still canonically definable and can be used for coding just like any real. More generally, small inner models seem mostly to present structural regularity almost as strong as that of $L$. Although time hasn't permitted, it would be an interesting line of research whether the study of gaps in some weaker or more local form is possible in $L[U]$, which doesn't satisfy Condensation or $GCH$ everywhere.

The inner models most relevant to the whole enterprise of Set Theory are canonical, and these present such regularities in some form even when big. So we are faced with a trade-off between foundationally relevant canonicity on the one hand, and analytic richness and intuitive behaviour on the other.

Regarding sets of reals, it is not clear at all that the existence of a Sierpi\'{n}ski set or a Luzin set should be a relevant argument against $CH$. Conversely, it is very desirable that non-measurable and non-Baire sets have high complexity, and we need very complex inner model constructions to ensure that.

The model theoretic trade-off we are forced to face is a common situation in foundational mathematics: we strive for a necessarily imperfect balance between the formal rigor of mathematical systems and the content of our intuitive concepts. The history of Logic in the last century demonstrates that only the dialogue between these two forces, the use of both, can push human mathematics forward.

\newpage\pagestyle{empty}
%%%%%%%%%%%%%%%%%%%%%%%%%%%%%%%%%%%%%%%%%%%%%%%%%%%%%%%%%%%%%%%%%%%%%%%%%
\backmatter
\selectlanguage{english}
\addcontentsline{toc}{chapter}{Bibliography}
\bibliographystyle{amsplain}

\end{document}